\documentclass[11pt]{article}
\usepackage{amssymb}    
\usepackage{amssymb,amsmath,amsthm,amscd}
\allowdisplaybreaks 

\begin{document}           
\newtheorem{proposition}{Proposition}[section]
\newtheorem{theorem}[proposition]{Theorem}
\newtheorem{definition}{Definition}[section]
\newtheorem{corollary}[proposition]{Corollary}
\newtheorem{Lemma}[proposition]{Lemma}
\newtheorem{conjecture}{Conjecture}

\newtheorem{Remark}{Remark}[section]
\newtheorem{example}{Example}[section]

\newcommand{\ds}{\displaystyle}
\renewcommand{\ni}{\noindent}
\newcommand{\pa}{\partial}
\newcommand{\Om}{\Omega}
\newcommand{\om}{\omega}
\newcommand{\va}{\varepsilon}
\newcommand{\var}{\varphi_{y,\la}}

\newcommand{\la}{\lambda}
\newcommand{\sik}{\sum_{i=1}^k}
\newcommand{\vov}{\Vert\omega\Vert}
\newcommand{\Umy}{U_{\mu_i,y^i}}
\newcommand{\lamns}{\lambda_n^{^{\scriptstyle\sigma}}}
\newcommand{\chiomn}{\chi_{_{\Omega_n}}}
\newcommand{\ullim}{\underline{\lim}}
\newcommand{\bsy}{\boldsymbol}
\newcommand{\mvb}{\mathversion{bold}}
\newcommand{\R}{{\mathbb R}}
\newcommand{\bR}{{\mathbb R}}
\newcommand{\bC}{\mathbb{C}}
\newcommand{\bE}{\mathbb{E}}
\newcommand{\bH}{\mathbb{H}}
\newcommand{\bP}{\mathbb{P}}
\newcommand{\cF}{\mathcal{F}}

\newcommand{\beq}{\begin{eqnarray*}}
\newcommand{\eeq}{\end{eqnarray*}}

\newcommand{\ben}{\begin{enumerate}}
\newcommand{\een}{\end{enumerate}}

\newcommand{\beqs}{\begin{eqnarray*}&\displaystyle}
\newcommand{\eeqs}{&\end{eqnarray*}}

\renewcommand{\theequation}{\thesection.\arabic{equation}}
\catcode`@=11 \@addtoreset{equation}{section} \catcode`@=12



\begin{center}
{\LARGE \bf Boundary Harnack inequalities for  regional fractional
Laplacian\footnote{Supported by NSFC/10501048, EPSRC/T26368/01,
EPSRC/D05379X/1}}
\end{center}
\begin{center} { \bf Qing-Yang Guan
  \footnote{ Institute of Applied Mathematics, Academy of Mathematics and Systems  Science, Chinese Academy of Sciences,   email: guanqy@amt.ac.cn
                      }}
\end{center}

\

\begin{abstract}
 We consider     boundary Harnack inequalities for
  regional fractional Laplacian which  are generators
 of  censored  stable-like processes on
 $G$ taking $$\kappa(x,y)/|x-y|^{n+\alpha}dxdy,\ \ \ x,y\in G$$ as the jumping measure.
  When $G$ is a $C^{1,\beta-1}$ open set, $1<\alpha<\beta\leq 2$
and  $\kappa\in C^{1}(\overline{G}\times \overline{G})$    bounded
between two positive numbers, we prove a  boundary Harnack
inequality giving $dist(x,\partial G)^{\alpha-1}$  order decay for
harmonic functions near the boundary. For  a $C^{1,\beta-1}$ open
set $D\subset \overline{D}\subset G $, $0<\alpha\leq (1\vee
\alpha)<\beta\leq 2$, we prove a boundary Harnack inequality giving
$dist(x,\partial D)^{\alpha/2}$ order decay for harmonic functions
near the boundary. These inequalities  are generalizations of  the
known results for the homogeneous case on $C^{1,1}$ open sets. We
also prove the boundary Harnack inequality for   regional fractional
Laplacian on Lipschitz domain.
\medskip

  { \bf MR(2000)} Subject Classification: Primary  31A20,  Secondary
 60G52, 60J75
\end{abstract}

\section{\normalsize Introduction }
 Let     $G$ be an open set   in $\R^n$ and $\kappa$   a
positive symmetric function on $\overline{G}\times \overline{G}$.
For $0< \alpha<2$,
  the regional fractional (fractional-like) Laplacian is defined  by
\begin{align} \label{oo}\Delta^{\frac{\alpha}{2},\kappa}_{\overline{G}}
u(x)=\lim_{\varepsilon\downarrow0}\mathcal{A}(n,-\alpha)\int_{y\in
G,|y-x|>\varepsilon}\frac{\kappa(x,y)(u(y)-u(x))}{|x-y|^{n+\alpha}}\
dy,\ \ \ x\in\overline{G},
\end{align}provided the limit  exists, see    \cite{GQY1}. Here
$\ds\mathcal{A}(n,-\alpha)=
 {|\alpha|2^{\alpha-1}\Gamma( ({n+\alpha})/{2})}{\pi^{- {n}/{2}}/\Gamma(1-
 {\alpha}/{2})}$
coming  from  $\Delta^{\frac{\alpha}{2},\kappa}_{\R^n}=-(-\Delta)^{
{\alpha}/{2}}$
 when $\kappa\equiv 1$ and we refer to Guan and Ma \cite{GQY} for $\kappa\equiv 1$ in (\ref{oo}). Under some regularity conditions, it is known that the
  $\alpha/2$ power of a second order elliptic operator with Neumann boundary condition is an example    of
(\ref{oo}). Since the integral kernel in (\ref{oo}) may not be
homogeneous in space, these operators to fractional Laplacian are
 similar to  the second order elliptic operators to     Laplacian. For
$1<\alpha<2$, among others, an explicit boundary Harnack inequality
(BHI) for $\Delta^{\frac{\alpha}{2},1}_{
 {\overline{G}}}$
was proved   in  Bogdan,   Burdzy and    Chen \cite{BBC} on
$C^{1,1}$ open sets, where it  is called the BHI of  the censored
stable processes. The main aim of this paper is to consider the same
type inequality for the nonhomogeneous case and the corresponding
  BHI on Lipschitz domain.

Boundary Harnack inequalities  are important tools in studying the
boundary value problems in partial differential equations and
   potential theory of Markov processes. Analytically, such
inequalities  describe  an uniform asymptotic behavior for solutions
of the
  Dirichlet problems near the
boundary.  In Chen and Kim \cite{CHKA},   the BHI in \cite{BBC} was
used in  the proof of  the Green function estimates of censored
stable processes. See also Bogdan \cite{BOG1} for  the Brownian
motion case. We refer to Bass \cite{BAS}, Chen, Kim and Song
\cite{CPS1} \cite{CPS2} for more applications of the BHI.

Boundary Harnack inequalities were first proved for Laplacian  in
Dahlberg \cite{DL}   and  Ancona \cite{A} on Lipschitz domains. It
was later extended to the  second order elliptic operators in
divergence form in Caffarelli,  Fabes, Mortola and  Salsa
\cite{CFMS}, and in nondivergence form  in Fabes, Garofalo,
Mar\'in-Malave and Salsa \cite{F}. A probability method for studying
such inequalities began in Bass and  Burdzy \cite{BB}.  This method
was applied to prove
 the BHI for Laplacian   on H\"older domains
for elliptic operators in divergence form in   Banuelos,   Bass
  and   Burdzy \cite{BBB}.

The study of  the BHI for the fractional Laplacian   began in Bogdan
\cite{BOGA}, Bogdan and Byczkowski \cite{BOG} on Lipschitz open
sets. Significant progresses have been made  on open sets in Song
and Wu \cite{SOG} and the recent paper Bogdan, Kulczycki and
Kwasnicki \cite{KTM}. An explicit BHI for the fractional Laplacian
was first given  in Chen and Song \cite{CS} on $C^{1,1}$ open sets.
Due to the jumps of stable processes or equivalently the nonlocal
property of their generators, the corresponding harmonic functions
show different feathers  from the Laplacian case.  Compared with the
fractional Laplacian case in \cite{CS} ($0<\alpha<2$), the BHI in
\cite{BBC} for the regional fractional Laplacian ($1<\alpha<2$)
gives  a different decay for harmonic functions near the boundary,
i.e., the former is of order $\rho(x)^{\alpha/2}$ and the later is
of order $\rho(x)^{\alpha-1}$.  In \cite{BBC}, the
    Markov processes associated with the  regional fractional Laplacian  under
   the  Dirichlet boundary condition were  first introduced and called  the  censored
stable processes. We refer to   \cite{BBC} for some other
characterizations of this process.

A standard box method to prove the BHI includes comparison of
harmonic measures around boxes, the Harnack inequality and the
Carleson estimate.  We refer to Bass and Burdzy \cite{BB}
  for this method in the diffusion case. The proof of the BHI in \cite{BBC} studied   these steps  mainly  by
  explicit
harmonic functions given in the same  paper and a relation between
the censored stable processes and the symmetric $\alpha$-stable
processes.   Some strong techniques are involved in this original
proof. Due to the importance of this inequality, it is helpful  to
  simplify  the  proof in \cite{BBC} and to study this
result in more general situations. In particular, our arguments can
be used to study the Lipschitz case which is an open problem in this
direction.  We remark that the (super sub) harmonic functions given
in \cite{BBC} plays a fundamental role in this paper.

To further    introduce the  results and the methods of this paper,
we prepare some definitions below.
 For $x=(x_1,\cdots,x_n)$ $\in \R^n$, we write $x=(\widetilde{x},x_n)$. Let $0<\gamma\leq 1$ and
$\Gamma:\R^{n-1}\rightarrow\R$. We say that $\Gamma$ is a
$C^{1,\gamma}$ function if it is differentiable and
\begin{align}\label{dsaf}
\|\Gamma\|_{1,\gamma}:=\sup_{\widetilde{y}\neq\widetilde{x},|\widetilde{y}-\widetilde{x}|
\leq 2}\frac{|\nabla\Gamma(\widetilde{y})-
\nabla\Gamma(\widetilde{x})|}{|\widetilde{y}-\widetilde{x}|^\gamma}<\infty,\end{align}
where $\nabla=(\partial/\partial x_i)_{i=1}^{n-1}$. The constant 2
in (\ref{dsaf}) is only for the convenience of the later use.  Let
$G$ be an open set in $\R^n$. We say that $G$ is a special
$C^{1,\gamma}$ domain if for some $C^{1,\gamma}$ function
$\Gamma:\R^{n-1}\rightarrow\R$, $G$ can be represented  as
$\{x=(\widetilde{x},x_n)\in \R^n,x_n>\Gamma(\widetilde{x})\}$. In
this case $G$ is also denoted by $G_\Gamma$. We say that $G$ is
$C^{1,\beta-1}$ if there exist $r_0>0$ and $\Lambda>0$ such that for
each $z\in
\partial
 G$, we can find a $C^{1,\beta-1}$ function $\Gamma_z:\R^{n-1}\rightarrow
 \R$ with $\|\Gamma\|_{1,\gamma}\leq\Lambda$ and an
 orthonormal coordinate system $CS_z$ such that
 \begin{align}\label{sorry}
G\cap B \bigr(z,r_0\bigr)=\{y=(y_1,\cdots,y_n):\ y_n> \Gamma_z(
y_1,\cdots y_{n-1}) \ \}\cap B\bigr( z,r_0\bigr).
\end{align}
 By rotation
and translation,  we can  always  assume that $\nabla \Gamma_z
(\widetilde{z})= \Gamma_z (\widetilde{z})=0$ in $CS_z$. The pair
$(r_0,\Lambda)$ is called the characteristics of $G$. The
characteristics of a Lipschitz open set is defined in a similar way.
For each $\delta>0$, set
\begin{align}
\ G_{\delta}'=\{y\in G: \rho(y)<\delta\},\ \ \ \ G_{\delta}=\{y\in
G: \rho(y)>\delta\},\end{align}
 where
$\rho(y)=dist(y,\partial G)$.

 Let $1<\alpha<2$ and   $\psi_1,\psi_2$ be positive
  functions in $ C^1(\overline{G}\times \overline{G})$.      Let $\kappa$ be a symmetric function
  on
$ {\overline{G}}\times \overline{{G}}$ taking values between two
positive numbers $C_1$ and $C_2$. Assume also that  for some
constant $C'>0$ and $\delta\in(0,r_0) $
\begin{align}\label{1}  \left\{
\begin{array}{r@{\quad \quad}l}  \left| \kappa(x,y)-\psi_1(x,y)-\psi_2(x,y)\frac{|x-y|^{n+\alpha}}{|x-\overline{y}|^{n+\alpha}}
  \right|
  \leq C'|x-y|\ ,\ \ \ \ \ \ \ \ \ \ \ x,y\in  G_{\delta}',\\
  {|\kappa(x,y)-\kappa(x,x)|}\leq C'{ |x-y|},\ \ \ \ \
\ \ \ \ \  \ \ \ \ \ \ \ \ \ \  \ \ \ \ \ \ \ \ \ \ \ \ \ \   \
x,y\in
 {G_{\delta/2}} ,
        \end{array} \right.\end{align}
       where $\overline{y}$ is the reflection point of $y$ with respect to
$\partial G$ (see section 4). Write
$$M:=C'+\sup_{x,y\in
\overline{G}, |x-y|<1
}(|\nabla_y\psi_1(x,y)|+|\nabla_y\psi_2(x,y)|)$$ and denote by
$(X_t)$ the reflected  stable-like process. The following theorem
 is an extension of the BHI in \cite{BBC} for $\kappa\equiv1$ on
$C^{1,1}$ domain.
\begin{theorem}\label{pr2.3.2} Assume that  $\alpha,\kappa$ satisfy all the conditions  above
and $G$ is  a $C^{1,1}$ open set in $\R^n$ with characteristics
$r_0\leq1$ and $\Lambda$.  Let $Q\in\partial G$ and $r\in (0,r_0)$.
Assume that
 $u\geq0 $ is a   function on $G$ which is not identical to
zero, harmonic on $G\cap B(Q,r)$ and vanishes continuously on
$\partial G\cap B(Q,r)$. Then there is a constant
$C=C(n,\alpha,\Lambda,\delta,C_1,C_2,M)$ such that
\begin{align}\label{yesr3} \frac{u(x)}{u(y)}\leq
C\frac{\rho(x)^{\alpha-1}}{\rho(y)^{\alpha-1}},\ \ \  \ \ x,y\in
G\cap B(Q,r/2).
\end{align}
Moreover, if $\psi_2\equiv0$ in (\ref{yesr3}), this boundary Harnack
inequality holds for  $C^{1,\beta-1}$ open sets with $1<\alpha<\beta
\leq 2$.
\end{theorem}
\noindent Here the notation
$C=C(n,\alpha,\Lambda,\delta_1,C_1,C_2,M)$ means that the constant
$C$ is   positive and depends only on   parameters in the bracket.
This convention will be used throughout the paper. When
$\psi_2\equiv0$ in ({\ref{yesr3}}), the last conclusion in Theorem
1.1 was conjectured
  in \cite{BBC}.
  We remark that  $\beta=\alpha$ is the critical value in our proof and
  Theorem 1.1 may not hold for this value.

In \cite{BBC}, when $ G$ is a special $C^{1,1}$ domain and
$\kappa\equiv1$, to get sharp estimates of harmonic measures,
(super) subharmonic functions are constructed by explicit harmonic
functions on $\R^n_+$ and non-explicit perturbations
  defined through the symmetric $\alpha$-stable process on
$\R^n$.     Here  we       construct
  explicit     (super) subharmonic functions by perturbation directly. This
  construction may  also be   used to prove the known  explicit BHI for Laplacian,
  i.e., $\alpha=\beta=2$ in (\ref{yesr3}).

For    the  Harnack inequality, we  may      adopt the method in
Bass and  Levin \cite{BAL}. Here we give another proof which might
be more straightforward for these nonlocal operators.  This proof is
similar to the proof of
  the Carleson estimate given in Lemma \ref{car} which is   an
  application of the
  box method  for  jump processes taking  (\ref{wsww}) as the key observation . We remark
that the method in \cite{BAL} can be applied  to prove the Harnack
inequality for jump diffusions, see Song and  Vondracek \cite{SV}.
Therefore by applying the method in this paper, we may prove the BHI
for  operators like $\Delta+\Delta^{\alpha/2}$ on $C^{1,1}$ open
sets, where  the decay is of order $\rho(x)$ near the boundary.

The BHI for the fractional  Laplacian on $C^{1,1}$ open sets was
proved by Poisson kernel estimates in \cite{CS}.  This and many
other  estimates  of the symmetric stable processes given  before
depend on their  explicit Poisson kernel and Green function  for  a
unit ball which are not available for the nonhomogeneous case. In
Lemma \ref{this},   we present explicit (super, sub) harmonic
functions of the fractional Laplacian on half spaces which  allow us
to study the nonhomogeneous case. See Theorem \ref{guan}. We notice
that the harmonic function in Lemma \ref{this} has been obtained in
 Banuelos and Bogdan \cite{BAB}. As applications, for the
fractional-like Laplacian under condition (\ref{dprf}), we may get
the sharp estimates of their Green function and  Poisson kernel as
in \cite{CS}    and hence we may get the BHI in \cite{KTM} under the
same conditions.

Another main result of this paper is the following boundary Harnack
inequality on Lipschitz domain. The strategy of the proof is
essentially the same as the proof of Theorem 1.1.
\begin{theorem}\label{pr2.3sadf.2} Let $G$ be a Lipschitz open
set in $\R^n$ with characteristics $r_0\leq1$ and $\Lambda$. Assume
that  $1<\alpha<2$ and $\kappa$ be a $C^1(\overline{G}\times
\overline{G})$ function bounded between two positive numbers. Let
$Q\in\partial G$ and $r\in (0,r_0)$.
  Then there is a constant
$C $ such that
\begin{align} \frac{u(x)}{u(y)}\leq
C\frac{v(x) }{v(y) },\ \ \  \ \ x,y\in G\cap B(Q,r/2),
\end{align}
where $u,v\geq0 $  are   functions on $G$ which is not identical to
zero, harmonic on $G\cap B(Q,r)$ and vanishes continuously on
$\partial G\cap B(Q,r)$.
\end{theorem}

To prove Theorem 1.2,   the heat kernel estimate of the reflected
stable processes in Chen and
  Kumagai \cite{CHKB}  is used to give some hitting probability
estimate. The censored stable processes can be extended to the
reflected
  processes  on $\overline{G}$ which is formulated  in \cite{BBC}.  For   general $\kappa$, these two
processes are called the
   censored stable-like process and the reflected
stable-like process respectively (see Remark 2.4 \cite{BBC}). They
are symmetric Markov processes on $G$ and $\overline{G}$,
respectively. The Dirichlet form of the reflected stable-like
process is
\begin{align}
\mathcal{E}^\kappa\bigr(u,v\bigr)&=
\frac{1}{2}\mathcal{A}(n,-\alpha)\int\int_{\overline{G}\times
\overline{G} }
\frac{\kappa(x,y)(u(x)-u(y))(v(x)-v(y))}{|x-y|^{n+\alpha}}\ dxdy,\nonumber \\
\mathcal{F}^\kappa&=  \{u\in L^2(\overline{G}):\ \ \
\mathcal{E}^\kappa\bigr(u,u\bigr)<\infty \},
\end{align}
where     $\kappa(x,y)$ is bounded between two   positive numbers
and $G$ is Lipschitz.  In   \cite{CHKB},   this reflected  process
was refined to be a Feller process $(X_t)$ on $\overline{G}$ under a
more general condition.

 The structure of this paper is the following. In Section 2 we study (super, sub) harmonic functions. In Section
 3 we prove  the   Harnack inequality. In Sections 4 and 6 we prove    boundary Harnack inequalities
 for the regional fractional-like Laplacian and the fractional-like Laplacian, respectively.
   The Lipschitz case is studied in section 5.    For $a,b\in\R$, $a\vee b:=\max\{a,b\}$.
We use $m(\cdot)$ to denote  the area measure of $(n-1)$-dimensional
subset. For any   set $U$, denote $\tau_{U}=\inf\{t>0: X_t\notin
U\}$.    For function $u$ on $\R^n$, we always take it as a function
on  $G$ by restriction when considering $\Delta_{G}^{ {\alpha}/{2}}
u$. The dimension $n$ is assumed bigger than two throughout the
paper.

\section{\normalsize    (Super) subharmonic functions for regional fractional Laplacian, $\kappa\equiv1$  }
 Let $u$ be a Borel function on $G$ and let $U$
be an open subset of $G$. We say that $u$ is a (super, sub) harmonic
function on $U$ with respect to $\Delta^{\frac{\alpha}{2},\kappa}_G$
if $\Delta^{\frac{\alpha}{2},\kappa}_Gu(x)(\leq,\geq)=0$ for $x\in
U$. We say that $u$ is a (super, sub) harmonic function on $U$ with
respect to the reflected stable-like process $(X_t)$ if
\begin{align} \label{hc}u(x)(\geq,\leq)= {E}_xu(X_{\tau_B}),\ \
 \  x\in B
\end{align}for any bounded open set $B$ with $\overline{B}\subset
U$, where $\tau_B=\inf\{t>0:X_t\notin B\}$. Under the conditions of
$\kappa$ and $G$ in this paper, the harmonic function for the
reflected stable-like process is continuous on $U$ (see Corollary
\ref{le5.1.8new} below). This implies that it is   harmonic for
$\Delta^{\frac{\alpha}{2},\kappa}_G$ in the weak sense (cf. Theorem
6.6  \cite{GMA} for $\kappa\equiv1$). When $u$ is $C^2$ and $\kappa$
is $C^1$,  Theorem 4.8 in \cite{GQY1} shows that these two
definitions are equivalent.  In what follows  we     use (\ref{hc})
for the definition of harmonic functions.   We write
 $\Delta^{\frac{\alpha}{2},1}_{
 {\overline{G}}}$
  by $\Delta^{{\alpha}/{2}}_{
 {G}}$ when $x\in G$.

In \cite{BBC}, to establish the BHI for $\kappa\equiv1$, the
following estimates are given for the regional fractional Laplacian
acting on function $u=\rho^{\alpha-1}$:
\begin{align}\label{2.1}\left\{
\begin{array}{r@{\quad \quad}l} \Delta_{G}^{ {\alpha}/{2}} u
(x)\leq A
\rho(x)^{\beta-2},\  \ \   x\in G_{1/A}', \ \ \ if\ \beta<2,\\
|\Delta_{G}^{ {\alpha}/{2}} u (x)|\leq A |\ln\rho(x)|,\  \
  x\in G_{1/A}', \ \ \ if\ \beta=2,
        \end{array}                  \right. \end{align}
where  $G$ is a special $C^{1,\beta-1}$ domain and  $A$ is a
positive constant ($1<\alpha< \beta< 2$). When $\beta<2$, we can not
get a bound for $\Delta_{G}^{ {\alpha}/{2}} u$ because  it  may take
$-\infty$. This is related  to the fact   that
 $\rho$ may not  be   $C^1$ when $\partial G$
is  $C^{1,\beta-1}$ ($\beta<2$). To improve   estimates (\ref{2.1}),
we replace the distance function by a ``height function" which is
equal to $x_n-\Gamma_z(x) $ in a neighborhood of $z\in\partial G$.

  Since  there is no difference for
   boundary conditions   discussed here when $n=1$, we
always assume that  $n\geq 2$. Let $w_p(y)=y_n^{p}$ for $
y\in\R^n_+$ and $p\in \R$.  Our starting point is the following
explicit harmonic functions and (\ref{1r12}) below given in  (5.4)
\cite{BBC}\begin{align}\label{111} &\Delta_{\R^n_+}^{
{\alpha}/{2}}w_{\alpha-1}(x)=0,\ \ \   x\in \R^n_+,\ \ \alpha\in
(1,2).
\end{align}
In the following lemma, when necessary, a function defined on a
domain is also   considered as a function on $\R^n$ by taking zero
outside.

\begin{Lemma}\label{hpr2.3.2}
Let $1<\alpha<\beta\leq2$ and let $\Gamma:\R^{n-1}\rightarrow\R$ be
a $C^{1,\beta-1}$ function with $\Gamma(\widetilde{0})=0$ and
$\nabla\Gamma(\widetilde{0})=0$. Define function
$h_{\alpha-1}(x)=(x_n-\Gamma(\widetilde{x}))^{\alpha-1}I_{\{|\widetilde{x}|<2\}}$
for $x\in D:=D_\Gamma$. Then there exists constant
$A_1=A_1(n,\alpha,\beta,\|\Gamma\|_{1,\beta-1})$ such that
\begin{align}\label{no}
|\Delta_{D}^{ {\alpha}/{2}} h_{\alpha-1} (x)|\leq\left\{
\begin{array}{r@{\quad \quad}l}  A_1
\rho(x)^{\beta-2},\  \ \   x\in D_{1}', \ |\widetilde{x}| <1,\ \ \ if\ \beta<2,\\
A_1 (|\ln\rho(x)|+1),\  \ \ \   x\in D_{1}', \ |\widetilde{x}| <1,\
\ \ if\ \beta=2.
        \end{array}                  \right. \end{align}
\end{Lemma} \noindent{\bf Proof} Denote $h_{\alpha-1}$ by $h$.
\ We only prove the lemma for $\alpha<\beta<2$ because  the proof
for $\beta=2$ is similar. Let $x\in D_1'$ with $ |\widetilde{x}| <1$
and choose a point $x_0\in\partial D$ satisfying
$\widetilde{x}=\widetilde{x_0}$. Denote by $\overrightarrow{n}(x_0)$
the inward   unit normal  vector at $x_0$ for $\partial D$ and set
$\Phi(y)=\langle y-x_0,\overrightarrow{n}(x_0)\rangle$ for $y\in
\R^n$. It is clear that $\Pi=\{y:\Phi(y)=0\}$ is the plane
 which is tangent to $\partial D$ at point $x_0$.
 Let
  $\Gamma^*:\widetilde{x}\in\R^{n-1}\rightarrow \R$ be the function of  plane
  $\Pi$, i.e.,
  $$
\langle(\widetilde{x},\Gamma^*(\widetilde{x}))-x_0,\overrightarrow{n}(x_0)\rangle=0,
  $$ and set
\begin{align} U=\{y=(\widetilde{y},y_n): y\in D,\ \
|\widetilde{y}-\widetilde{x}|<1,\ \
y_n<2+2^{\beta}\|\Gamma\|_{1,\beta-1}\ \}. \ \ \nonumber
 \end{align}
Write  $\overline{h}(y)=|y_n-\Gamma^*(\widetilde{y})|$ for $y\in
\R^n$. Applying the assumption that  $\partial G$ is $C^{1,\beta-1}$
and $\nabla\Gamma(\widetilde{x})-\nabla\Gamma^*(\widetilde{x})=0$,
we have by the mean value theorem
\begin{align}\label{hhhh}
|\overline{h}({y})-h^{\frac{1}{\alpha-1}}({y})|\leq|\Gamma(\widetilde{y})-\Gamma^*(\widetilde{y})|\leq
\|\Gamma\|_{1,\beta-1}|\widetilde{y}-\widetilde{x}|^{\beta},\ \ \ \
\ y\in U.
\end{align}
Let $\rho_{_{\Pi}}(y)=dist(y,\Pi)$  for $y\in \R^n$ and
${D}_{\Gamma^*}=\{y\in \R^n:y_n>\Gamma^*(\widetilde{y})\}$. It is
clear that $\overline{h}={\sqrt{1+|\nabla
\Gamma({\widetilde{x_0}})|^2}\rho_{_{\Pi}}}$. So we have by
(\ref{111})
\begin{align}\label{forgqq}
\Delta^{
{\alpha}/{2}}_{D_{\Gamma^*}}\overline{h}^{\alpha-1}(y)={({1+|\nabla
\Gamma({\widetilde{x_0}})|^2})^{\frac{\alpha-1}{2}}}\Delta^{\frac{\alpha}{2}}_{D_{\Gamma^*}}\rho_{_{\Pi}}^{\alpha-1}(y)=0,\
\ \   y\in {D}_{\Gamma^*}.
\end{align}
Denote $$A=\{y: \Gamma^*(\widetilde{y})<y_n<\Gamma(\widetilde{y}),\
|\widetilde{y}-\widetilde{x}|<1\}\cup\{y:
\Gamma(\widetilde{y})<y_n<\Gamma^*(\widetilde{y}),\
|\widetilde{y}-\widetilde{x}|<1\}.$$ Noticing that
$\overline{h}^{\alpha-1}(x)=h(x)$ and $B(x,1)\cap D \subset U$ (by
the fact that  $x_n\leq 1+2^{\beta}\|\Gamma\|_{1,\beta-1}$), we have
by (\ref{forgqq})\begin{align}
\label{shan}&\lim_{\varepsilon\downarrow 0}\bigg|\int_{y\in
D,|y-x|>\varepsilon}\frac{h(y)-h(x)}{|x-y|^{n+\alpha}}\
dy\bigg|\nonumber\\\leq
&\limsup_{\varepsilon\downarrow0}\bigg|\int_{y\in
U,|y-x|>\varepsilon}\frac{\overline{h}^{\alpha-1}(y)-\overline{h}^{\alpha-1}(x)}{|x-y|^{n+\alpha}}\
dy\bigg|+\limsup_{\varepsilon\downarrow0}\bigg|\int_{y\in
U,|y-x|>\varepsilon}\frac{h(y)-\overline{h}^{\alpha-1}(y)}{|x-y|^{n+\alpha}}\
dy\bigg|\nonumber\\+&\limsup_{\varepsilon\downarrow0}\bigg|\int_{y\in
D\setminus U,|y-x|>\varepsilon}\frac{h(y)-h(x)}{|x-y|^{n+\alpha}}\
dy\bigg|\nonumber\\\leq&\int_{A}\frac{|\overline{h}^{\alpha-1}(y)-\overline{h}^{\alpha-1}(x)|}{|x-y|^{n+\alpha}}\
dy+\int_{B(x,1)^c}\frac{|\overline{h}^{\alpha-1}(y)-\overline{h}^{\alpha-1}(x)|}{|x-y|^{n+\alpha}}\
dy\nonumber\\+&\int_{U}\frac{|h(y)-\overline{h}^{\alpha-1}(y)|}{|x-y|^{n+\alpha}}\
dy+\int_{B(x,1)^c}\frac{ |h(y)-h(x)|}{|x-y|^{n+\alpha}}\ dy
\nonumber\\:=&I_1+I_2+I_3+I_4.
\end{align}
Noticing that  $A\subset \{y: |y_n-(x_0)_n|\leq
2^{\beta-1}\|\Gamma\|_{1,\beta-1}|\widetilde{y}-\widetilde{x_0}|\
\}$, we have
\begin{align}
|x-y|\geq(1+2^{2\beta-2}\|\Gamma\|^2_{1,\beta-1})^{-1/2}\overline{h}(x)
\geq(1+2^{\beta-1}\|\Gamma\|_{1,\beta-1})^{-1}\overline{h}(x),\ \
  y\in A,\nonumber
\end{align}which implies
$$|x-y|\geq\frac{(1+2^{\beta-1}\|\Gamma\|_{1,\beta-1})^{-1}{\overline{h}(x)}+|\widetilde{y}-\widetilde{x}|}{2},
\ \ \ y\in A.$$    By (\ref{hhhh}), we also have
$\overline{h}(y)\leq
\|\Gamma\|_{1,\beta-1}|\widetilde{y}-\widetilde{x}|^\beta$ for $y\in
A$. Therefore \begin{align} I_1\leq&\int_{ 0}^{1}\
dr\int_{|\widetilde{y}-\widetilde{x}|=r}I_A(y)
\frac{|\overline{h}^{\alpha-1}(y)-\overline{h}^{\alpha-1}(x)|}{|x-y|^{n+\alpha}}
\ m(d{y})\nonumber\\\leq&(\|\Gamma\|_{1,\beta-1}^{\alpha-1}+1)\int_{
0}^{\overline{h}(x)^{\frac{1}{\beta}}\wedge1}dr\int_{|
\widetilde{y}-\widetilde{x}|=r}I_A(y)\frac{\overline{h}(x)^{\alpha-1}}{|x-y|^{n+\alpha}}\
m(dy)\nonumber\\+&
(\|\Gamma\|_{1,\beta-1}^{\alpha-1}+1)\int^{1}_{\overline{h}(x)^{\frac{1}{\beta}}\wedge1}
dr\int_{|\widetilde{y}-\widetilde{x}|=r}I_A(y)\frac{{r^{\beta(\alpha-1)}}}{|x-y|^{n+\alpha}}\
m(dy)\nonumber\\\leq&(2\pi)^n(\|\Gamma\|_{1,\beta-1}^{\alpha}+\|\Gamma\|_{1,\beta-1})\int_{
0}^{\overline{h}(x)^{\frac{1}{\beta}}\wedge1}{\overline{h}(x)^{\alpha-1}}{(\frac{r+(1+2^{\beta-1}\|\Gamma\|_{1,\beta-1})^{-1}
{\overline{h}(x)}}{2})^{-\alpha+\beta-2}}\ dr\nonumber\\+& (2\pi)^n
(\|\Gamma\|_{1,\beta-1}^{\alpha}+\|\Gamma\|_{1,\beta-1})\int^{1}_{\overline{h}(x)^{\frac{1}{\beta}}\wedge1}
{r^{\beta(\alpha-1)}}{(\frac{r+(1+2^{\beta-1}\|\Gamma\|_{1,\beta-1})^{-1}{\overline{h}(x)}}{2})^{-\alpha+\beta-2}}\ dr\nonumber\\
\leq&(2\pi)^n(\|\Gamma\|_{1,\beta-1}^{\alpha}+\|\Gamma\|_{1,\beta-1})2^{\alpha-\beta+2}
\frac{(1+2^{\beta-1}\|\Gamma\|_{1,\beta-1})^{\alpha-\beta+1}}{\alpha-\beta+1}\overline{h}(x)^{\beta-2}\nonumber\\
+&(2\pi)^n(\|\Gamma\|_{1,\beta-1}^{\alpha}+\|\Gamma\|_{1,\beta-1})2^{\alpha-\beta+2}\int^{1}_{\overline{h}(x)^{\frac{1}{\beta}}\wedge1}
r^{\alpha\beta-\alpha-2} dr.
\end{align}
As ${\alpha\beta-\alpha-2}>\beta^2-2{\beta}-1$ for
$1<\alpha,\beta\leq 2$, we get
\begin{align}\label{2.8}
\int^{1}_{\overline{h}(x)^{\frac{1}{\beta}}\wedge1}
r^{\alpha\beta-\alpha-2} dr\leq
\int^{1}_{\overline{h}(x)^{\frac{1}{\beta}}\wedge1}
r^{\beta^2-2\beta-1} dr \leq \frac{1}{2\beta-\beta^2}
(\overline{h}(x)\wedge1)^{\beta-2}.
\end{align}The following properties follows from the  definitions of $\overline{h}$   and
$h$.
\begin{align}\label{name}
|\overline{h}(x)-\overline{h}(y)|&\leq(1+\|\Gamma\|_{1,\beta-1})|x-y|,\
\ \ \ \  \ \ \ \ \ \ \ \ \ \ \ \ y\in \R^n,\\ \label{one}
|h(x)^{\frac{1}{\alpha-1}}-h(y)^{\frac{1}{\alpha-1}}|&\leq(1+2^{\beta-1}\|\Gamma\|_{1,\beta-1})|x-y|,\
\ \  \ \ \ \ \ \ \   y\in D,\ |\widetilde{y}|\leq 2,\\\label{OOO}
\rho(y)^{\alpha-1}\leq h(y)&\leq
(1+2^{\beta-1}\|\Gamma\|_{1,\beta-1})^{\alpha-1}\rho(y)^{\alpha-1},\
\ \ \ y\in D_1',\ |\widetilde{y}|\leq 1.
\end{align}Noticing  that $\rho(x)<1$ and $h(y)=0$ for $|\widetilde{y}|>2$, by (\ref{name})-(\ref{OOO}) \begin{align}\label{2.12}
&I_2+I_4\nonumber\\
\leq&\int_{B(x,1)^c}\frac{(1+\|\Gamma\|_{1,\beta-1})^{\alpha
-1}}{|x-y|^{n+1}}\
dy+\int_{B(x,1)^c}\frac{(1+2^{\beta-1}\|\Gamma\|_{1,\beta-1})^{\alpha-1}}{|x-y|^{n+1}}\
dy+\int_{B(x,1)^c}\frac{h(x)}{|x-y|^{n+\alpha}}\ dy\nonumber\\\leq&
(2\pi)^n\big((1+\|\Gamma\|_{1,\beta-1})^{\alpha-1}+(1+2^{\beta-1}\|\Gamma\|_{1,\beta-1})^{\alpha-1}
+ {(1+2^{\beta-1}\|\Gamma\|_{1,\beta-1})^{\alpha-1}}\big). \
\end{align}
To estimate $I_3$ we define a transform
$\Psi(y)=(\widetilde{z},z_n)$ by
$$\widetilde{z}=\widetilde{y},\ \ \ z_n=y_n-\Gamma^*(\widetilde{y}),\ \ \ \
y\in \R^n.$$ We see that $|\frac{\partial \Psi}{\partial y}|=1$,
where $\frac{\partial \Psi}{\partial y}$ is the Jacobian determinant
of $\Psi$. We can also check  that
$$|y_1-y_2|\geq(1+\|\Gamma\|_{1,\beta-1})^{-1}
|\Psi(y_1)-\Psi(y_2)|,\ \ \ \ \ for\ y_1,y_2\in \R^n,$$ and
$\Psi(U)\subset \{y:|\widetilde{y}-\widetilde{x}|<1, \ |y_n|\leq
2^{\beta+1}(\|\Gamma\|_{1,\beta-1}+1) \}$. Hence by  (\ref{hhhh}),
the inequality
\begin{align}\label{2.13}
|b^{\alpha-1}-a^{\alpha-1}|\leq b^{\alpha-2}|b-a|,\ \ \
b>0,a>0,1<\alpha<2,
\end{align}
and applying the transform $\Psi$,  we have
 \begin{align}\label{2.14}
I_3\leq&\int_{
U}\frac{(1+\|\Gamma\|_{1,\beta-1})^{n+\alpha+1}|\overline{h}(y)|^{\alpha-2}
|\widetilde{y}-\widetilde{x}|^{\beta}}{|\Psi(x)-\Psi(y)|^{n+\alpha}}\
dy \nonumber\\=&\int_{\Psi(U)}\
\frac{(1+\|\Gamma\|_{1,\beta-1})^{n+\alpha+1}|z_n|^
{\alpha-2}|\widetilde{z}-\widetilde{x}|^{\beta}}{|(\widetilde{x},\overline{h}(x))-z|^{n+\alpha}}\
dz\nonumber\\\leq
&\int_{-2^{\beta+1}(\|\Gamma\|_{1,\beta-1}+1)}^{2^{\beta+1}(\|\Gamma\|_{1,\beta-1}+1)}\
dr\int_{z_n=r,|\widetilde{z}-\widetilde{x}|<|r-\overline{h}(x)|}
\frac{(1+\|\Gamma\|_{1,\beta-1})^{n+\alpha+1}|r|^{\alpha-2}|\widetilde{z}
-\widetilde{x}|^{\beta}}{|(\widetilde{x},\overline{h}(x))-z|^{n+\alpha}}\
m(dz)\nonumber\\+&
\int_{-2^{\beta+1}(\|\Gamma\|_{1,\beta-1}+1)}^{2^{\beta+1}(\|\Gamma\|_{1,\beta-1}+1)}\
dr\int_{z_n=r,|r-\overline{h}(x)|\leq|\widetilde{z}-\widetilde{x}|\leq1}
\frac{(1+\|\Gamma\|_{1,\beta-1})^{n+\alpha+1}|r|^{\alpha-2}}{|\widetilde{z}-\widetilde{x}|^{n+\alpha-\beta}}\
m(dz)\nonumber\\
\leq&\frac{(2\pi)^n}{\beta+1}\int_{0}^{2^{\beta+1}(\|\Gamma\|_{1,\beta-1}+1)}\
\frac{(1+\|\Gamma\|_{1,\beta-1})^{n+\alpha+1}|r|^{\alpha-2}}{|\overline{h}(x)-r|^{\alpha-\beta+1}}\
dr\nonumber\\
+&\frac{(2\pi)^n}{\alpha-\beta+1}\int_{0}^{2^{\beta+1}(\|\Gamma\|_{1,\beta-1}+1)}\
\frac{(1+\|\Gamma\|_{1,\beta-1})^{n+\alpha+1}|r|^{\alpha-2}}{|\overline{h}(x)-r|^{\alpha-\beta+1}}\
dr\nonumber\\
\leq&\frac{2(2\pi)^n}{\alpha-\beta+1}\bigg(\int_0^{2\overline{h}(x)}\frac{(1+\|\Gamma\|_{1,\beta-1})^{n+\alpha+1}}{r^{2-\alpha}|r-\overline{h}(x)|^{\alpha-\beta+1}}\
dr+ \int_{2\overline{h}(x)}^{2^{\beta+1}
(\|\Gamma\|_{1,\beta-1}+1)}\frac{(1+\|\Gamma\|_{1,\beta-1})^{n+\alpha+1}}{(r-\overline{h}(x))^{-\beta+3}}\
dr\bigg)\nonumber\\
\leq&\frac{2(2\pi)^n(1+\|\Gamma\|_{1,\beta-1})^{n+\alpha+1}}{\alpha-\beta+1}\bigg(\int_0^{2}\frac{\overline{h}(x)^{\beta-2}}{r^{2-\alpha}|r-1|^{\alpha-\beta+1}}\
dr+\frac{1}{2-\beta}\overline{h}(x)^{\beta-2}\bigg).
\end{align}
  Combining (\ref{shan})-(\ref{2.8}),(\ref{2.12}) and (\ref{2.14}), we get (\ref{no}).
\medskip\qed

\begin{Remark}\label{RE}\  Estimates
 (\ref{no}) may  not hold if we take $\beta=\alpha$ in Lemma \ref{hpr2.3.2}. For $n=2$,
 $\Gamma(x_1)=|x_1|^{\beta}$ and $x^*=(0,t)$ with  $t>0$, we can check that
$\int_{U}\frac{h(y)-\overline{h}^{\alpha-1}(y)}{|x^*-y|^{2+\alpha}}\
dy=-\infty$ and $I_1,I_2,I_4$ are all finite. This gives
$\Delta_{D_\Gamma}^{ {\alpha}/{2}} h (x^*)=- \infty$.  When
$\alpha<\beta< 2$,  we can also prove that
$\Delta_{D_\Gamma}^{\frac{\alpha}{2}} \rho^{\alpha-1}(x^*) $ may
take $-\infty$. We still consider the above example. Let $x^*_0$ be
the point on $\partial D_\Gamma$ such that $|x_0^*-x^*|=\rho(x)$ and
$(x_0^*)_1>0$. Let $$ U=\{(y_1,y_2):\ \ y_2> |y_1|^\beta\ or\
y_1\leq 0\}\cap \{(y_1,y_2):y_2>0\}$$ and denote the distance
function to $\partial U$ by $\xi(x)$.  Since  $\xi$ is smooth in a
neighborhood of $x^*$, we know  that $\Delta_{U}^{\frac{\alpha}{2}}
\xi^{\alpha-1} (x^*)$ is   finite. On the other hand,
$$\int_{D_{\Gamma}}\frac{\rho(y)^{\alpha-1}-\xi(y)^{\alpha-1}}{|x^*-y|^{2+\alpha}}dy
=-\infty.$$Hence we have  $\Delta_{D_\Gamma}^{\frac{\alpha}{2}}
\rho^{\alpha-1}(x^*) =-\infty$.
\end{Remark}

Recall
 $w_p(y)=y_n^{p}$ for $ y\in\R^n_+$.
By (5.4) in \cite{BBC}\begin{align}\label{1r12} &\Delta_{\R^n_+}^{
{\alpha}/{2}}w_p(x)=\mathcal{A}(n,-\alpha)
\frac{\omega_{n-1}}{2}\mathcal{B}(\frac{\alpha+1}{2},\frac{n-1}{2})\gamma(\alpha,p)x^{p-\alpha},\
\ \   x\in \R^n_+,\ \ p\in (-1,\alpha),
\end{align}
where $\omega_{n-1}$ is the $(n-2)$-dimensional Lebesgue measure of
the unit sphere in $\R^{n-1}$,  $\mathcal{B}$ is the Beta function
and
$\gamma(\alpha,p)=\int_0^1\frac{(t^p-1)(1-t^{\alpha-p-1})}{(1-t)^{1+\alpha}}\
dt$. In what follows we denote the constant on the right hand side
of (\ref{1r12}) by $C(n,\alpha,p)$.
\begin{Lemma}\label{pr2.3.2y} Let $\alpha$, $\Gamma$ and $D$ be
described in  Lemma \ref{hpr2.3.2} and let $p$ be a number such that
$\alpha> p>\alpha-1$. Define function
$h_p(x)=(x_n-\Gamma(\widetilde{x}))^{p}I_{\{|\widetilde{x}|<2\}}$ on
$D=D_\Gamma$. Then there exists constant
$A_2=A_2(n,\alpha,\beta,p,\|\Gamma\|_{1,\beta-1})$ such that
\begin{align}\label{211} \Delta_{D}^{ {\alpha}/{2}} h_p (x)
\geq A_2 \rho(x)^{p-\alpha},\  \ \ \   x\in D_{1/A_2}', \
|\widetilde{x}| <1.\end{align}
\end{Lemma} \noindent{\bf Proof}
\  We use the definitions and  the notations in the proof of Lemma
\ref{hpr2.3.2}. Following the arguments in  (\ref{shan}), for $x\in
D_1'$ with $ |\widetilde{x}| <1$ we have
\begin{align}
\label{shanxu}&\lim_{\varepsilon\downarrow0}\int_{y\in
D,|y-x|>\varepsilon}\frac{h_p(y)-h_p(x)}{|x-y|^{n+\alpha}}\
dy\nonumber\\\geq&C(n,\alpha,p)x^{p-\alpha}-\int_{A}\frac{|\overline{h}^{p
}(y)-\overline{h}^{p }(x)|}{|x-y|^{n+\alpha}}\
dy-\int_{B(x,1)^c}\frac{|\overline{h}^{p }(y)-\overline{h}^{p
}(x)|}{|x-y|^{n+\alpha}}\ dy\nonumber\\-&\int_{
U}\frac{|h_p(y)-\overline{h}^{p }(y)|}{|x-y|^{n+\alpha}}\
dy-\int_{B(x,1)^c}\frac{|h_p(y)-h_p(x)|}{|x-y|^{n+\alpha}}\
dy\nonumber\\=&C(n,\alpha,p)x^{p-\alpha}-I_1-I_2-I_3-I_4.
\end{align}By   similar calculations  as in Lemma \ref{hpr2.3.2}, we can find
constant $k_1$ such that
\begin{align}\label{shanxu1}
&I_1\leq k_1(\rho(x)^{\beta+p -\alpha-1}\vee1+|\ln \rho(x)|),\nonumber\\
\ &I_3 \leq k_1 (\rho(x)^{\beta+p -\alpha-1} \vee 1+|\ln \rho(x)|),
\ \ \ \ I_2+I_4\leq k_1.
\end{align}
Noticing  that $p-\alpha<0$, $\beta>1$ and $C(n,\alpha,p)>0$, we
obtain (\ref{211}) by (\ref{shanxu}) and
(\ref{shanxu1}).\medskip\qed

\begin{Lemma}\label{pr}
Let $\alpha$, $\Gamma$, $D$, $h_{\alpha-1}$ and $h_p$ be objects
described in Lemmas \ref{hpr2.3.2} and \ref{pr2.3.2y}. Let $f$ be a
bounded function in $ C^1(\overline{D})$. Then there exists
constant $\ds A_3=A_3(n,\alpha,p,\sup_{y\in \R^n}|
f(y)|,\sup_{|y|<2}|\nabla f(y)|)$ such that for $x\in D\cap B(0,1)$
\begin{align}\label{no212} \int_D\frac{|(f(y)-f(x))(h_p(y)-h_p(x))|}
{|y-x|^{n+\alpha}}\ dy\leq    \left\{
\begin{array}{r@{\quad \quad}l}  A_3 (|\log\rho(x)|+1),\  \ \ \  p=\alpha-1,\\
A_3,\  \ \ \ \ \ \ \ \ \ \ \ \ \ \ \ \  \alpha-1<p<\alpha .
        \end{array} \right.
            \end{align}

\end{Lemma} \noindent{\bf Proof} We only prove the lemma  for
$p=\alpha-1$ because the others  can be proved  similarly. Denote
$h_{\alpha-1}$ by $h$ and let $x\in D\cap B(0,1)$. By
(\ref{one}),(\ref{OOO}) and (\ref{2.13}),
\begin{align}  &\int_D\frac{|(f(y)-f(x))(h(y)-h(x))|}{|y-x|^{n+\alpha}}\ dy
\nonumber\\
\leq&\sup_{|y|<2}|\nabla f(y)|\int_{D\cap
\{\rho(x)<|y-x|\leq1\}}\frac{|h(y)^{\frac{1}{\alpha-1}}-h(x)^{\frac{1}{\alpha-1}}|^{\alpha-1}
}{|y-x|^{n+\alpha-1}}\ dy\nonumber\\+&\sup_{|y|<2}|\nabla
f(y)|\int_{D\cap
B(x,\rho(x))}\frac{|h(y)^{\frac{1}{\alpha-1}}-h(x)^{\frac{1}{\alpha-1}}|h(x)^{\frac{\alpha-2}{\alpha-1}}
}{|y-x|^{n+\alpha-1}}\ dy\nonumber\\+& 2\sup_{y\in \R^n}|
f(y)|\int_{D\cap B(x,1)^c\cap\{y:|\widetilde{y}|\leq2\}}
\frac{|h(y)^{\frac{1}{\alpha-1}}-h(x)^{\frac{1}{\alpha-1}}|^{\alpha-1}
}{|y-x|^{n+\alpha}}\ dy\nonumber\\+&2\sup_{y\in \R^n}|
f(y)|\int_{\{y:|\widetilde{y}|>2\}}\frac{h(x)}{|y-x|^{n+\alpha}}\
dy\nonumber\\
\leq&\sup_{|y|<2}|\nabla f(y)|\int_{D\cap
\{\rho(x)<|y-x|\leq1\}}\frac{(1+2^\beta\|\Gamma_z\|_{1,\beta-1})^{\alpha-1}}{|y-x|^{n}}\
dy\nonumber\\+& \sup_{|y|<2}|\nabla f(y)|\int_{D\cap
B(x,\rho(x))}\frac{(1+2^\beta\|\Gamma_z\|_{1,\beta-1})^{\alpha-1}\rho(x)^{\alpha-2}}{|y-x|^{n+\alpha-2}}\
dy \nonumber\\+&2\sup_{y\in \R^n}| f(y)|\int_{D\cap
B(x,1)^c}\frac{(1+2^\beta\|\Gamma_z\|_{1,\beta-1})^{\alpha-1}}{|y-x|^{n+1}}\
dy+2\sup_{y\in \R^n}|
f(y)|\frac{(1+2^\beta\|\Gamma_z\|_{1,\beta-1})^{\alpha-1}}{\alpha}(2\pi)^n\nonumber\\
\leq &(\sup_{|y|<2}| \nabla f(y)|+2\sup_{y\in \R^n}|
f(y)|){(1+2^\beta\|\Gamma_z\|_{1,\beta-1})^{\alpha-1}}(2\pi)^n( -\ln
\rho(x)+\frac{1}{2-\alpha}+1+\frac{1}{{\alpha}}),\nonumber
\end{align} which completes the proof.\qed\medskip

For  $C^1$  function $\kappa$ on $\overline{G}\times \overline{G}$
and $Q \in \partial G$, in the proposition below we denote
\begin{align}\label{C_0}C_0=\sup_{x,y\in B(Q,1)\cap G }|\nabla_y\kappa(x,y)|.\end{align}
\begin{proposition}\label{pr2.4.} Let $1<\alpha<\beta\leq 2$ and $G$
  a $C^{1,\beta-1}$ open set in $\R^n$ with characteristics $r_0=1$
and $\Lambda$.  Let $\kappa$ be a   $C^1$ function    on
$\overline{G}\times \overline{G}$  taking values between two
positive numbers $C_1$ and $C_2$. Then for $\alpha-1\leq p< \alpha$
and $Q\in
\partial G$, there
exist  function $u_p$ and positive constants $A_4=A_4( \Lambda)$,
$A_5=A_5\big(n,\alpha,\beta,p,\Lambda,C_0,C_1,C_2 \big)$ such that
\begin{align}\label{r}A_4^{-1}I_{G\cap B(Q,2/3)}\rho(x)^p\leq u_p(x)\leq A_4I_{G\cap
B(Q,2/3)}\rho(x)^p,\ \ \ x\in G,\end{align} and
\begin{align}\label{no212} \Delta_{G}^{ {\frac{\alpha}{2}},\kappa} u_p (x)
\geq& A_5 \rho(x)^{p-\alpha},\  \ \ \  \ \ \ \   x\in G\cap B(Q,1/A_5),\ \ \alpha-1< p< \alpha,       \\
|\Delta_{G}^{  {\frac{\alpha}{2}},\kappa} u_{\alpha-1}
(x)|\leq&\left\{
\begin{array}{r@{\quad \quad}l}  A_5
\rho(x)^{\beta-2},\  \ \   x\in G\cap B(Q,1/2),\ \ \ if\ \beta<2,\\
A_5 |\ln\rho(x)|,\  \ \ \   x\in G\cap B(Q,1/2),\ \ \ if\ \beta=2.
        \end{array} \right.
            \end{align}
        \end{proposition}
        \noindent{\bf Proof}\  Without loss of generality, we assume that $Q=0$ and
        take the
        coordinate system  $CS_Q$ (see (\ref{sorry})). Define functions
         $u_{p}(x)=(x_n-\Gamma_Q(\widetilde{x}))^{p}I_{{G\cap\{B(Q,2/3)\}}}$
         on $G$
        and $ v_{p}(x)=(x_n-\Gamma_Q(\widetilde{x}))^{p}I_{{  |\widetilde{x}|<2
        }}$ on $D_{\Gamma_Q}$
         for $\alpha-1\leq p<\alpha$. It is easy to see that  (\ref{r}) holds. When $\kappa \equiv1$,  noticing that for $x\in G\cap B(0,1/2)$ the
        integral in (\ref{oo}) for $u_p$ on $G\cap B(Q,2/3)^c$ and $v_p$ on $D_{\Gamma_Q}\cap B(Q,2/3)^c$ can be   bounded
        by constants depending on $n$ and $ \alpha$, we can prove this proposition  by Lemma
        2.1 and Lemma \ref{pr2.3.2y}.
       For    general cases, the conclusion   can be proved by the case $\kappa \equiv1$, Lemma 2.3
       and the following identity:
\begin{align}\label{dd}
&\Delta_{G}^{\frac{\alpha}{2},\kappa} h (x)\nonumber\\
=& \mathcal{A}(n,-\alpha)\lim_{\varepsilon \downarrow 0}\int_{y\in
G, |y-x|>\varepsilon}
\frac{(\kappa(x,y)-\kappa(x,x))(h(y)-h(x))}{|x-y|^{n+\alpha}}dy
+\kappa(x,x)\Delta_{G}^{ {\alpha}/{2}} h (x).
\end{align}\qed

\section{\normalsize   Harnack inequalities of $\Delta^{\frac{\alpha}{2},\kappa}_G$ }

  The following   example can be found
in \cite{GQY1}.
\begin{example}Let  $G=\R^n_+$ and  $\overline{y}=(\widetilde{y},-y_n)$ for
$y=(\widetilde{y},y_n)$.  For
$\kappa(x,y)=1+\frac{|x-y|^{n+\alpha}}{|x-\overline{y}|^{n+\alpha}}$,
$\Delta^{\frac{\alpha}{2},\kappa}_{G}$ is the formal generator of
the subordinate  reflected Brownian motion on $\overline{\R^n_+}$.
When $G=(0,1)$ and $\kappa(x,y)=\sum_{k=-\infty }^{\infty}
{|x-y|^{1+\alpha}}/{|x\pm y+2k|^{1+\alpha}}$,
$\Delta^{\frac{\alpha}{2},\kappa}_{G}$ is the formal  generator  of
the subordinate reflected Brownian motion on $[0,1]$.
\end{example}
\begin{Remark}
Define function
 $w_p(y)=y_n^{p}$ for $ y\in\R^n_+$.
When
$\kappa(x,y)=\frac{|x-y|^{n+\alpha}}{|x-\overline{y}|^{n+\alpha}}$,
we have  (see\cite{GQY1})
\begin{align}\label{112}
&\Delta_{\R^n_+}^{\frac{\alpha}{2},\kappa}w_p(x)=\mathcal{A}(n,-\alpha)
\frac{\omega_{n-1}}{2}\mathcal{B}(\frac{\alpha+1}{2},\frac{n-1}{2})\overline{\gamma}(\alpha,p)x^{p-\alpha},\
\ \   x\in \R^n_+,\ \ p\in (-1,\alpha),
\end{align}
where $\overline{\gamma}(\alpha,p)=\int_0^1
{(t^p-1)(1-t^{\alpha-p-1})}/{(1+t)^{1+\alpha}}\ dt$. This gives the
same (super,sub) harmonic functions as the homogeneous case in
(\ref{1r12}), which will be used later.
\end{Remark}

Notice that the derivatives of $\kappa$ in the   examples above  are
not bounded. To give results including   these examples we introduce
the following condition.
 Let $0<C_1<C_2$, $C_3>0$ and
$\gamma\leq 0$. We say that  $\kappa$ or  the   reflected
 stable-like process $(X_t)$ satisfies condition $[C_1,C_2,C_3,
\gamma]$ if \begin{align}\label{dsf} C_1<\kappa(x,y)<C_2,\ \ x,y\in
\overline{G};\ \ \ \
|\kappa(x,y)-\kappa(x,x)|<C_3(\rho(x)^\gamma\vee1)|x-y|,\   \ x,y\in
G.\end{align} We can check  that   functions $\kappa$ in the Example
3.1 above satisfy  condition
  $[C_1,C_2,C_3,-1]$ for some constants $C_1,C_2,C_3>0$.

 Next  we prepare a stochastic calculus
formula for $(X_t) $.  For a measurable function $f$ on $G$, denote
$f\in\mathcal{L}^1_u(G)$ if \begin{align}\label{rrr} \sup_{x\in G}
\int_{G} \frac{|f(x)-f(y)|}{(1+|x-y|)^{n+\alpha}}\ dy<\infty.
\end{align}For any   subset $U\subseteq \R^n$ and $0<\gamma\leq 1$, we say
that $u$ is uniformly $\gamma$-H\"{o}lder continuous on $U$ if
\begin{align}\label{rr}\sup_{(y,z)\in U\times
U}\frac{|u(y)-u(z)|}{|z-y|^{\gamma}}<\infty.\end{align}  We shall
denote  ($u\in C^{1,\gamma}(U)$) $u\in C^{\gamma}(U)$ if (all the
first derivatives of $u$) $u$ is uniformly $\gamma$-H\"{o}lder
continuous on $U$. For any $\delta>0$ and $A\subseteq\R^n$, define
$\tau_A=\inf\{t>0: X_t\in A^c\}$ and $A^\delta=\{y: |y-x|<\delta,\
for \ some \ x\in A\}$. For any relatively open subset  $A$ of
$\overline{G}$, we denote by $(p_t^{A}) $ and $G^{A} $ the
probability  transition function  and the Green function    of
$(X_t) $  killed upon leaving ${A}$, respectively. In  \cite{GQY1},
a semi-martingale decomposition of $f(X_t)$ is given for $f\in
C^{2}_c(\overline{G})$ (see  \cite{GQY} for the homogeneous case).
To consider more general functions, we prove the following results.

\begin{proposition}\label{pr4.1.3}
Let  G be   a Lipschitz open set in
 $\R^n$. For  $1\leq \alpha<\beta\leq 2$, let  $\kappa$ be  a symmetric  function on
  ${\overline{G}}\times{\overline{G}}$
 satisfying condition $[C_1,C_2,C_3, \gamma]$  with $0\geq\gamma> \alpha-3$. For
  $0<\alpha< \beta\leq 1$,  let  $\kappa$ be  a measurable
symmetric function on   ${\overline{G}}\times{\overline{G}}$ bounded
between $C_1$ and $C_2$.  Then for
 $f$ belonging to
 \begin{align}\label{rp}
C^{1,\beta-1}(\overline{G})\cap\mathcal{L}^1_u(G), \ \ \ 1\leq
\alpha<\beta\leq 2;\ \ \ \
C^{\beta}(\overline{G})\cap\mathcal{L}^1_u(G), \ \ \ 0<
\alpha<\beta\leq 1,\end{align} we have
\begin{align}\label{formula}
&f(X_t)=f(x_0)+M_t+\int_0^t
\Delta^{\frac{\alpha}{2},\kappa}_{G}f(X_s)\ ds, \ \ \ a.s.\ \
  x_0\in \overline{G},
\end{align}
where   $(M_t)_{t\geq 0}$ is a    martingale. If $A$ is a relatively
open set in $\overline{G}$ and, for some  $\delta>0$, $f$ satisfies
(\ref{rp}) with $ \overline{G}$ replaced by $ \overline{G} \cap
A^\delta $, then
\begin{align}\label{3.1fg}
E_{x_0}(f(X_{t\wedge\tau_A}))&=f(x_0)+ E_{x_0}( \int_0^{t\wedge
\tau_A}\Delta^{ \frac{\alpha}{2},\kappa}_Gf(X_s)ds),\ \ \  \ x_0\in
A, \ t\geq 0.\end{align}  Moreover, if $P_{x_0}(X_{\tau_A}\in
\partial A)=0$ and $f$ is a positive function such that $f=0$ on
$A$, then (\ref{3.1fg}) still holds.
\end{proposition}
\noindent{\bf Proof} Assume that $f$ satisfies (\ref{rp}). For
$1\leq \alpha<\beta\leq 2$, by (\ref{rr}), the derivatives of $f$ at
point $x\in G$ can be bounded by $k_1(1+|x|)$ for some constant
$k_1>0$. By this estimate, (\ref{rrr}) and straightforward
calculations for the integral in (\ref{oo}) on sets $B(x,\rho(x))$,
$ (G\cap B(x,1))\setminus B(x,\rho(x))$ and $G\setminus B(x,1)$
respectively, we can prove  that for some constant $k_2$,
\begin{align}\label{cc}|\Delta^{ \frac{\alpha}{2},\kappa}_Gf(x)|\leq k_2(1+|x|)
  \rho(x)^{(2+\gamma-\alpha)\wedge (1-\alpha)} ,\ \ \ \ \ x\in
G_1',\ \ 1\leq \alpha<\beta\leq 2.\end{align} By (\ref{rrr}) and
(\ref{rr}) we can   find constant $k_3$ such that
\begin{align}\label{c}|\Delta^{ \frac{\alpha}{2},\kappa}_Gf(x)|\leq
k_3,\ \ \ \ \ x\in  G_1,\ \ \ 1\leq \alpha<\beta\leq 2.\end{align}
By calculating   the integral in (\ref{oo}) on sets $G\cap B(x,1)$
and $G\setminus B(x,1)$ respectively,  we can also check
\begin{align}\label{cc,}  |\Delta^{ \frac{\alpha}{2},\kappa}_Gf(x)|\leq& k_4,\ \
\ \ \ \     x\in G,\ \  0< \alpha<\beta\leq 1,\end{align} for some
constant $k_4$. Noticing that $(2+\gamma-\alpha)\wedge
(1-\alpha)>-1$, with the help of the heat kernel estimates in
\cite{CHKB} and (\ref{cc})-(\ref{cc,}) we can prove that
$E_x(\int_0^1 |\Delta^{ \frac{\alpha}{2},\kappa}_Gf(X_t)|dt)$ is a
bounded function on $\overline{G}$ (c.f.   Lemma 4.6 \cite{GQY}).
This implies  that $E_x(\int_0^t |\Delta^{
\frac{\alpha}{2},\kappa}_Gf(X_t)|dt)$ is a bounded function on
$\overline{G}$   for any $t>0$. Thus we can prove (\ref{formula}) by
Theorem 5.25 \cite{FOT} at time  $t\wedge \tau_{B(0,n)}$ (c.f.
Theorem 4.1 \cite{GQY}) and letting $n\rightarrow \infty$. Formula
(\ref{3.1fg}) is a consequence of (\ref{formula}) by
 approximation procedure.
   \qed\medskip

For a  relatively  open set $A$ in $\overline{G}$, we say that $A$
has  outer cone property in $\overline{G}$ if, for some $\eta>0$ and
each $Q\in\partial A$, there is a cone in $\overline{G}\setminus A$
isometric to $\{x\in \R^n: |(x_1,\cdots,x_{n-1})|<\eta |x_n|\}$ and
taking $Q$ as the vertex.
\begin{proposition}\label{le5.1.8} Let $\alpha, G$ and $\kappa$ be
the same as in Proposition  \ref{pr4.1.3}.
 Let $A\subseteq\overline{G}$ be a  relatively  open set    with outer
cone property in $\overline{G}$ and define  $\tau=\inf\{t>0: X_t\in
{A}^c\}$. Then the distribution of $X_{\tau}$  is absolutely
continuous on $\overline{G}\setminus {A}$ when $(X_t) $ starting
from $A$. For any $t>0$, we have
\begin{align}\label{H2}
&P_x \{X_{\tau}I_{\{\tau\leq t\}}\in dy
\}/dy \nonumber\\
=&\mathcal{A}(n,-\alpha)\int_0^t \bigg(\int_{A }
\frac{\kappa(z,y)p^{A}(s,x,z)}{|z-y|^{n+\alpha}}dz\bigg)ds,\ \
  (x,y)\in A\times (\overline{G}\setminus  {A}).
\end{align}
 Furthermore,
 \begin{align}\label{H222}
&P_x \{X_{\tau} \in dy
\}/dy \nonumber\\
=&\mathcal{A}(n,-\alpha) \int_{A }
\frac{\kappa(z,y)G^{A}(x,z)}{|z-y|^{n+\alpha}}dz,\ \
  (x,y)\in A\times (\overline{G}\setminus  {A}).
\end{align}\end{proposition}\noindent{\bf Proof}
 $\ $ To show that $P_x\{X_{\tau}\in \partial A\}=0$ for
$x\in A$, by the method in Lemma 6 \cite{BOGA},   we only need to
prove that there exists a constant $c$ such that
$P_x\{\tau_{B(x,\rho(x))}\in G\setminus A \}>c$ for any $x\in {A}$
(the boundedness assumption in \cite{BOGA} is not necessary because
   $A$ can be approximated  by $A\cap B(0,n)$ by letting
$ n\uparrow\infty $). We omit the proof of this estimate because it
is similar to (\ref{H5}) below. Thus we can prove    (\ref{H2}) by
  Proposition \ref{pr4.1.3}.
Formula (\ref{H222}) is a consequence of (\ref{H2}).
 \qed\medskip

\begin{Lemma}\label{scaling} Let $0<\alpha<2$ and let G be a  Lipschitz open set in
 $\R^n$.  Assume that   $\kappa$  is   a symmetric  function on  ${\overline{G}}\times{\overline{G}}$
 satisfying condition $[C_1,C_2,C_3, -1]$.  Let $\lambda\geq 1$ and  define process $((Z_t)_{t\geq 0},
Q_{x_0})=((\lambda X_{\lambda^{-\alpha}t})_{t\geq
0},P_{x_0/\lambda})$ for $x_0\in \overline{\lambda G}$. Then $(Z_t)
$ is a reflected stable-like process on $\overline{\lambda G}$
satisfying condition $[C_1,C_2,C_3, -1]$.
\end{Lemma} \noindent{\bf Proof}
The conclusion   can be proved by checking   that the jumping
measure of $(Z_t) $ is $$ \ \ \ \ \ \ \ \ \ \ \ \ \ \ \ \ \ \ \ \ \
\ \ \ \ \ \ \ \ \ \ \ \ \ \ \ \ \ \ \ \ \ \ \ \ \ \ \
\frac{\kappa(x/\lambda,y/\lambda)}{|x-y|^{n+\alpha}}\ dxdy,\ \ \ \
x, y\in {\lambda G}.\ \ \ \ \ \ \ \ \ \ \ \ \ \ \ \ \ \ \ \ \ \ \ \
\ \ \ \ \ \ \ \ \ \ \ \ \ \ \ \ \ \ \ \ \ \ \ \ \qed$$

\begin{Lemma} \label{forget} Let $0<\alpha<2$ and let  G be a  Lipschitz open set in
 $\R^n$ with characteristics $r_0=1$ and $\Lambda$. Assume that    $\kappa$ is  a symmetric  function
 on  ${\overline{G}}\times{\overline{G}}$
 satisfying condition $[C_1,C_2,C_3, -1]$. Then for  $0<\varepsilon<1$,
  there exists constants $A_6=A_6(n,\alpha,  C_2,C_3,\varepsilon)$ and $A_6'=A_6'(n,\alpha, C_1 )$
   such that  for any $0<r\leq r_0/2$
\begin{align}\label{hitting time}
A_6r^{\alpha}\leq \inf_{y\in B(x,(1-\varepsilon)r)}E_{y}
\tau_{B(x,r)}\leq \sup_{y\in \overline{G}}E_{y} \tau_{B(x,r)}\leq
 A_6'r^{\alpha},\ \ \   x\in G\ with \ \rho(x)>2r.
\end{align}
Moveover, the last inequality in (\ref{hitting time}) holds   for
all  $x\in \overline{G}$ provided  $r<r_0/4$, where $A_6'$  depends
further  on $\Lambda$.
\end{Lemma}\noindent{\bf Proof} By the scaling property in Lemma
\ref{scaling} and the Lipschitz condition of $G$, we can assume that
$r=1$.
 Choose  $f_1\in C^2(\overline{G})$  such that $0\leq f_1\leq 1$
 and  $$f_1(y)=0,\ \ \  y\in B(x,1-\varepsilon);\ \ \ \ \ f_1(y)=1,\ \  \ y\in B(x,1-\varepsilon/2)^c.$$ By direct
  calculation, we can find a constant $k_1=k_1(n,\alpha, C_2,C_3,\varepsilon)$ such
 that $|\Delta_G^{\frac{\alpha}{2},\kappa}f_1(y)|<k_1$ for $y\in B(x,1)$. Thus we can prove
 the first  inequality in (\ref{hitting time}) by applying
  formula  (\ref{3.1fg}) to $E_y(f_1(X_{\tau_{B(x,1)}}))$. Similarly,
  with the help of Proposition \ref{le5.1.8}, the last inequality in  (\ref{hitting time})
  can be proved  by considering function $f_2=I_{\overline{G}\setminus B(x,1)}$, where
 we can check that $\Delta_G^{\frac{\alpha}{2},\kappa}f_2(y)>k_2$, $y\in \overline{G}\cap B(x,1)$, for some constant
 $k_2=k_2(n,\alpha, C_1 )$. For the last conclusion,
  $\overline{G}\cap B(x,1)$ may not have the outer cone property in $\overline{G}$, where
  we need to replace $B(x,1)$ by a bigger set in $B(x,2)$ satisfying  this property.     \qed
 \medskip

The next theorem extends the Harnack inequality for the   censored
stable process  in \cite{BBC}.
\begin{theorem}\label{pr4.1.33}Let $0<\alpha<2$ and let  G be a  Lipschitz open set in
 $\R^n$ with characteristics $r_0$ and $\Lambda$. Assume that    $\kappa$ is  a symmetric
  function  on  ${\overline{G}}\times{\overline{G}}$
 satisfying condition $[C_1,C_2,C_3, -1]$.
Let $0<r\leq 1$, $k\in \{1,2,\ldots\}$ and $x_1, x_2\in G$ such that
  $B(x_1,r)\cup B(x_2,r)\subset G$  and  $|x_1-x_2|<2^kr$. If $u\geq
0$ is harmonic for $(X_{t}) $ on $B(x_1,r)\cup B(x_2,r)$, then there
exists constant $A_{7}=A_{7}(n,\alpha,C_1,C_2,C_3)$ such that
\begin{align}\label{Harnack}
A_{7}^{-1}2^{-k(n+\alpha)}u(x_2)\leq u(x_1)\leq A_{7}
2^{k(n+\alpha)}u(x_2).
\end{align}
\end{theorem}\noindent{\bf Proof}$\ $ For simplicity we assume
$\kappa\equiv1$. Let  $y \in G$ with $\rho(y)\geq r$.  First we
prove that there exists a constant $k_1=k_1(n,\alpha)$ such that
\begin{align}\label{Harnack1}
  u(y_1)\leq k_1 u(y_2),\ \ \   y_1,y_2 \in
B(y,r/2),
\end{align}provided $u\geq 0$ is
harmonic for $(X_t) $ on $B(y,r)$. To show this we only need to
prove that
\begin{align}\label{Harnack1n}
  u(y_1)\leq k_1 u(y_2),\ \ \   y_1,y_2 \in
B(y,r/2)\ and\ |y_1-y_2|>r/3.
\end{align}
Approximating  by functions $u_k:=E_x\big({(u\wedge
k)(X_{\tau_{B(y,r)}})}\big)$, we can assume that $u$ is bounded.  By
scaling we can also assume $r=1$.

 Let $ y_1,y_2 \in B(y,1/2)$ such
that $ |y_1-y_2|>1/3$. Suppose that  $u(y_1)>Mu(y_2)$ for some big
number $M$ and   we can construct a sequence of points $(x_k)\in
B(y_1,1/6)$ such that $x_0=x_1=y_1$ and
\begin{align}\label{true}
|u(x_{k})|\geq (1+\delta)^{k-1}Mu(y_2),\ \  \ \ |x_{k}-x_{k-1}|\leq
12^{-1}(k-1)^{-2} ,\ \ \ k\geq 1,  \end{align} for some $\delta>0$,
then the contradiction between   the boundedness of $u$ and
$\lim_{k\rightarrow \infty}u(x_k)=\infty$ leads to (\ref{Harnack1n})
(here we assume that $u(y_2)>0$ because $u(y_2)=0$  implies that
$u\equiv0$ on $G\cap B(y,1)$, c.f. (\ref{t2}) below).

 Suppose
that $(\ref{true})$ holds     for      $k= 1$  and some $\delta $,
  $M $ which will be fixed later. Setting
$B_k=B(x_k,24^{-1}k^{-2})$ and $\tau_k=\tau_{B_k}$ for $k\geq 1$,
we have by Proposition \ref{le5.1.8} and Lemma \ref{forget}
\begin{align}\label{t} &P_{x_k}\{X_{\tau_k} \in
G\setminus(2B_k)\}\nonumber\\=&
\mathcal{A}(n,-\alpha)\int_{B_k}\int_{y\in G\setminus (2B_k)}
\frac{ G^{B_k}(x_k,z)}{|z-y|^{n+\alpha}} dzdy\nonumber\\
\geq &
2^{-(n+\alpha)}\mathcal{A}(n,-\alpha)\int_{B_k}G^{B_k}(x_k,z)\int_{y\in
G\setminus (2B_k)}
\frac{1}{|x_k-y|^{n+\alpha}}dzdy\nonumber\\
= & 2^{-(n+\alpha)}\mathcal{A}(n,-\alpha)(E_{x_k}{\tau_k})\int_{y\in
G\setminus (2B_k)} \frac{1}{|x_k-y|^{n+\alpha}}dy
\nonumber\\
\geq &k_1\end{align} for some constant $k_1=k_1(n,\alpha)$.
Similarly, by setting $B_0=B(y_2, 1/6)$ and $\tau_0=\tau_{B_0}$, we
also have
\begin{align}\label{t2} &u(y_2)=E_{y_2}(u(X_{\tau_0})I_{X_{\tau_0} \in
G\setminus B_0}) \geq k_2(n,\alpha) \int_{y\in G\setminus B_0}
\frac{u(y)}{|y_2-y|^{n+\alpha}}dy.
\end{align}
By (\ref{hc}) and an  estimate of $P_{y_2}(\tau_0\in
B(y,1/2)\setminus (2B_0))$ similar to (\ref{t}), we can find $y_3\in
B(y,1/2)\setminus (2B_0)$ such that $u(y_3)\leq k_3 u(y_2)$ for some
constant $k_3=k_3(n,\alpha)$. Similar to (\ref{t2}), we have
\begin{align}\label{t4} &u(y_3)  \geq k_4(n,\alpha) \int_{y\in   B_0}
\frac{u(y)}{|y_3-y|^{n+\alpha}}dy.
\end{align}

Noticing that $|y-x_k|\geq \frac{1}{12 k^2}(|y-y_2|\vee|y-y_3|)$ for
$y\in G\setminus(2B_k)$, we have by Proposition \ref{le5.1.8} and
Lemma \ref{forget}
\begin{align}\label{t1} &E_{x_k}(u(X_{\tau_k})I_{X_{\tau_k} \in
G\setminus(2B_k)})\nonumber\\
\leq  &
2^{n+\alpha}\mathcal{A}(n,-\alpha)(E_{x_k}{\tau_k})\int_{y\in
G\setminus(2B_k)} \frac{u(y)}{|x_k-y|^{n+\alpha}}dy
\nonumber\\
\leq &2^{n+\alpha}(\frac{1}{24k^2})^{\alpha}A_6'
\mathcal{A}(n,-\alpha) \int_{y\in G\setminus(2B_k)}
\frac{u(y)}{|x_k-y|^{n+\alpha}}dy\nonumber\\
\leq &k_5(n,\alpha)k^{2n} \big(\int_{y\in G\setminus B_0}
\frac{u(y)}{|y_2-y|^{n+\alpha}}dy+ \int_{y\in B_0}
\frac{u(y)}{|y_3-y|^{n+\alpha}}dy\big).
\end{align}
By (\ref{t2})-(\ref{t1}) and $u(y_3)\leq k_3 u(y_2)$, we have
$E_{x_k}(u(X_{\tau_k})I_{X_{\tau_k} \in G\setminus(2B_k)})\leq
k_6(n,\alpha)k^{2n} u(y_2)$. Thus by (\ref{hc}), (\ref{true}) and
(\ref{t}) we have
\begin{align}\label{t_3}
(1+\delta)^{k-1}Mu(y_2)\leq u(x_k)\leq (1-k_1)\sup_{y\in
(2B_k)\setminus B_k}u(y)+k_6 k^{2n} u(y_2).
\end{align}
Now   choose  $\delta=k_1/2$ and    $K_0=K_0(n,\alpha)\geq 1$ such
that for any $M>1$
\begin{align}\label{t6}(1+\delta)^{m-1}M-k_6 m^{2n} \geq
\frac{1-k_1}{1-k_1/2}(1+\delta)^{m-1}M,\ \ \ \ m\geq K_0.
\end{align} If $x_k$ with
$k>K_0$ satisfies  (\ref{t_3}) for some $M>1$, then (\ref{t_3}) and
(\ref{t6}) show  that there exists $x_{k+1}$ satisfying (\ref{true})
for
   $k+1$. By  (\ref{t_3}), we can also choose $M=M(n,\alpha)>1$ big enough such
   that (\ref{t6}) holds for $1\leq k\leq K_0$. Therefore, we can finish
   the proof of (\ref{Harnack1}) by induction.

  Next we assume that
$2^kr>|x_1-x_2|> r$. We have for $x\in B(x_1,r/3)$
\begin{align}\label{Hr} \Delta^{\frac{\alpha}{2},\kappa}_G
I_{B(x_2,r/3)}(x)=&\mathcal{A}(n,-\alpha)\int_{y\in
B(x_2,r/3)}\frac{1}{|x-y|^{n+\alpha}}\ dy\nonumber\\
\geq&k_7(n,\alpha) r^{-\alpha} 2^{-k(n+\alpha)}.
\end{align}
 By (\ref{3.1fg}),(\ref{hitting time}) and (\ref{Hr}) we have
\begin{align}\label{H5} &P_{x_1}\{X_{\tau_{B(x_1,r/3)}}\in  B(x_2,
r/3)\}\nonumber\\
=&E_{x_1}(\int_0^{\tau_{B(x_1,r/3)}}
\Delta^{\frac{\alpha}{2},\kappa}_G
I_{B(x_2,r/3)}(X_t)dt) \nonumber\\
\geq&k_7A_6  3^{-\alpha} 2^{-k(n+\alpha)}.
\end{align}
By (\ref{Harnack1}) and (\ref{H5}),
\begin{align}\label{H1} u(x_1)=& E_{x_1}(u(X_{\tau_{B(x_1, r/3)}}))
\geq k_8(n,\alpha) 2^{-k(n+\alpha)} u(x_2),
\end{align}
which completes the proof. \qed\medskip

\begin{corollary}\label{le5.1.8new} Let $\alpha, G$ and  $\kappa$ be the same as in Theorem \ref{pr4.1.33}
and let  $u$ be a $(X_t) $ harmonic function in an open subset $D$
of $G$. Then  $u$ is continuous on $D$.
\end{corollary}\noindent{\bf Proof} \ Let $x\in D$ and $\rho_D(x)=\inf\{|x-y|: y\in \partial D\}$. By
Theorem \ref{pr4.1.33}£¬ we see that $u$ is bounded on
$B(x,2\rho_D(x)/3)$. Set  $\tau=\tau_{B(x,\rho_D(x)/3)}$. By the
strong Markov property, we have
\begin{align}\label{yes}
u(y)=E_y  u(X_{\tau})=E_y [ u(x_{\tau})I_{t\geq \tau}]+E_y
[u(X_t)I_{t<\tau}].
\end{align}
By the  continuity of the heat kernel  in \cite{CHKB},  we see
$E_\cdot [u(X_t) I_{\{t<\tau\}}]\in C(B(x,\rho_D(x)/3))$(c.f.
Proposition 3.6   \cite{GMA}). On the other hand,
\begin{align}
&\bigg| E_y[u(X_\tau)I_{\{t\geq \tau\}}]\bigg|\nonumber\\\leq
&\left(\sup_{z\in (B(x,2\rho_D(x)/3))}|u(z)|\right)P_y\{t\geq \tau,
X_\tau \in (B(x,2\rho_D(x)/3))\}+\bigg| E_y[u(X_\tau)I_{{t\geq
\tau}}I_{X_\tau\notin B(x, 2\rho_D(x)/3)}]\bigg|.\nonumber
\end{align}
Therefore, by (\ref{H2}), (\ref{yes}), Theorem \ref{pr4.1.33} and
the dominated convergence theorem, we  need only to check that
$P_y\{t\geq \tau\}$ converges
 to zero uniformly on  $B(x, \rho_D(x)/3)$ when $t\downarrow0$.
This  follows from   facts that $ P_y \{t\geq \tau\}=1-P_y\{t<
\tau\}\in C_b(B(x, \rho_D(x)/3)) $ and $
\lim\limits_{t\downarrow0}P_y\{t\geq \tau\}=0$ for $ y\in B(x,
\rho_D(x)/3). $ \qed\medskip

\section{\normalsize  Boundary Harnack inequality of
  $\Delta^{\frac{\alpha}{2},\kappa}_G$   on $C^{1,\beta-1}$
($C^{1,1}$) open sets}

Next we assume  that $0\in \partial G$ and choose the  coordinate
system $CS_0$. For $x\in \R^n$, $r>0$, let $\triangle(x,a,r)$ be the
box defined by
\begin{align}
\triangle(x,a,r)=\{y=(\widetilde{y},y_n)\in G:
0<y_n-\Gamma_0(\widetilde{y})<a, |\widetilde{y}-\widetilde{x}|<r\}.
\end{align}
The following result is a special case  of Theorem 1.1.
\begin{Lemma}
\label{pr2.3.20}Let $1<\alpha<\beta\leq2$ and let $G$ be a
$C^{1,\beta-1}$ open set with characteristics  $r_0=1$ and
$\Lambda$. Assume that $\kappa$ satisfies the  conditions in
Proposition \ref{pr2.4.}. Then there exist constants
$A_8=A_8(n,\alpha,\beta,\Lambda,C_0,C_1,C_2)<1/2$ and
$A_9=A_9(n,\alpha,\beta,\Lambda, C_0,C_1,C_2)$ such that
\begin{align}\label{3.1101}
 A_9^{-1} \rho(x)^{\alpha-1} \leq &P_x\{X_{\tau_{\triangle(0,A_8,A_8)}}\in
\triangle(0,2A_8,A_8)\}\nonumber\\\leq&
P_x\{X_{\tau_{\triangle(0,A_8,A_8)}}\in G\}\ \leq A_9
\rho(x)^{\alpha-1}
\end{align}   for $ x\in \triangle(0,A_8,A_8)$ with
$\widetilde{x}=0$.
\end{Lemma}\noindent{\bf Proof}\   We assume that $\kappa
\equiv1$ because the proof is the same for  the general case.
 Let $
p=(\alpha-1+((\alpha+\beta-2)\wedge1))/2$  and define \begin{align}
v_1(y)=&u_{\alpha-1}(y)+u_p(y)  ,\nonumber
\end{align}where  $u_{\alpha-1}$ and $u_p$ are functions defined in Proposition
\ref{pr2.4.}.
Since $p-\alpha>\beta-2$, by Proposition \ref{pr2.4.}, there exists
$k_1=k_1(n,\alpha,\beta,\Lambda) $ such that
$\triangle(0,2k_1,2k_1)\subseteq B(0,r_0)$ and
\begin{align}\label{3.121}
 \Delta_{G}^{  {\alpha}/{2}}v_1(y)\geq0,\  \ \ \ \ \ \   y\in
\triangle(0,k_1,k_1).\end{align}

 Let $\phi$ be a $C^2$
function on $\overline{G} $ such that
$$\phi(y)=|\widetilde{y}|=y_1^2+\ldots+y_{n-1}^2,\ \ \  \ |y|<1;\ \ \ \ \ 1\leq \phi(y)\leq 2,
\ \ \  \ |y|\geq1.$$ Define \begin{align}
v_2(y)=&u_{\alpha-1}(y)-u_p(y)/(2A_4^2)+12k_1^{-2}A_4^3\phi(y).\nonumber
\end{align}
By Lemma 3.4  \cite{BBC}, we have $ |\Delta_{G}^{
 {\alpha}/{2}}\phi(y)|\leq k_2(\rho(y)^{1-\alpha}\vee1)$, $y\in \triangle(0,k_1,k_1)$, for  some
 constant
$k_2=k_2(n,\alpha)$. Thus by $p-\alpha<1-\alpha$ and Proposition
\ref{pr2.4.}, there exists $m=m(n,\alpha,\beta,\Lambda)\leq k_1/2$
such that
\begin{align}\label{3.12}
 \Delta_{G}^{  {\alpha}/{2}}v_2(y)\leq0,\  \ \ \ \ \ \   y\in
\triangle(0,m,k_1).\end{align}  Since $v_2\geq 3A_4^3$ on
$G\setminus\triangle(0,\infty,k_1/2)$ and $v_2(y)\leq
A_4\rho(y)^{\alpha-1}$ for $y\in G\cap B(0,r_0)$ with
$\widetilde{y}=0$,   we have by applying (\ref{3.1fg}) and
(\ref{3.12})
\begin{align}\label{3.10}
P_x\{X_{\tau_{\triangle(0,m,k_1/2)}}\in
G\setminus\triangle(0,\infty,k_1/2)\}\leq 3^{-1}A_4^{-2}
\rho(x)^{\alpha-1}
\end{align}
for $x\in \triangle(0,m,k_1)\ with \ \widetilde{x}=0$.

 Noticing that
$\sup_{y\in G}|v_1(y)|\leq 2A_4 $ and $v_1(y)\geq
A_4^{-1}\rho(y)^{\alpha-1}$  for $y\in G\cap B(0,r_0)$ with
$\widetilde{y}=0$,  by   (\ref{3.1fg}) and (\ref{3.121}), we have
\begin{align}\label{3.110}
  P_x\{X_{\tau_{\triangle(0,m,k_1)}}\in
G\setminus  \triangle(0,m,k_1/2)\}\ \geq 2^{-1}A_4^{-2}
\rho(x)^{\alpha-1},\ \ \  x\in \triangle(0,m,k_1)\ and \
\widetilde{x}=0.
\end{align}
 Combing  (\ref{3.10}) and (\ref{3.110}), we have
\begin{align}\label{label}P_x\{X_{\tau_{\triangle(0,m,k_1/2)}}\in
\triangle(0,\infty,k_1/2)\setminus\triangle(0,m,k_1/2)\}\geq
6^{-1}A_4^{-2} \rho(x)^{\alpha-1}
\end{align}
for  $x\in \triangle(0,m,k_1)\ with \ \widetilde{x}=0$. By
(\ref{label}) and (\ref{H222}), we can find a  constant
$k_2=k_2(k_1,m,\Lambda)$
\begin{align}\label{labeltt}P_x\{X_{\tau_{\triangle(0,m,k_1/2)}}\in
\triangle(0,k_1,k_1/2)\setminus\triangle(0,m,k_1/2)\}\geq
k_2A_4^{-2} \rho(x)^{\alpha-1}.
\end{align}

By (\ref{labeltt}),    for $x\in \triangle(0,m,k_1)$ with
$\widetilde{x}=0$
\begin{align}&\label{aa} P_x\{X_{\tau_{\triangle(0,k_1,k_1)}}\in
\triangle(0,2k_1,k_1)\setminus\triangle(0,k_1,k_1)\}\nonumber\\\geq&
P_x\{X_{\tau_{\triangle(0,m,k_1/2)}}\in
\triangle(0,2k_1,k_1)\setminus\triangle(0,k_1,k_1)\}\nonumber\\\geq&
P_x\{X_{\tau_{\triangle(0,m,k_1/2)}}\in
\triangle(0,k_1,k_1/2)\setminus\triangle(0,m,k_1/2)\}\cdot\nonumber\\
& \sup_{y\in
\triangle(0,k_1,k_1/2)\setminus\triangle(0,m,k_1/2)}P_y\{X_{\tau_{\triangle(0,k_1,k_1)}}\in
\triangle(0,2k_1,k_1)\setminus\triangle(0,k_1,k_1)\}\nonumber\\
\geq &k_2k_3A_4^{-2}\rho(x)^{\alpha-1},
\end{align} where we use the fact that for some $k_3=k_3(k_1,m,\Lambda)$
\begin{align}\label{a} &P_y\{X_{\tau_{\triangle(0,k_1,k_1)\ }}\in
\triangle(0,2k_1,k_1)\ \setminus \triangle(0,k_1,k_1)\ \}\geq
k_3,\nonumber\\ & \   \ \ \ \ y\in \triangle(0,k_1,k_1/2)\ \setminus
\triangle(0,m,k_1/2).
\end{align}which can be proved by the same calculation as  (\ref{H5}).
Setting $A_8=k_1$,   (\ref{3.10}), (\ref{aa}) and (\ref{a}) yield
the first inequality of (\ref{3.1101}) for $x\in
\triangle(0,k_1,k_1)\ with \ \widetilde{x}=0$.

 Set $v_3(x)=v_2(x)I_{x\in G, |x|<1/2}+I_{x\in G, |x|\geq 1/2}$, by
 Proposition
\ref{pr2.4.} we can choose $k_1$ small enough such that
\begin{align}\label{3.12}
 \Delta_{G}^{  {\alpha}/{2}}v_3(y)\leq0,\  \ \ \ \ \ \   y\in
\triangle(0,k_1,k_1).\end{align} This estimate and Proposition
\ref{pr4.1.3} gives the second inequality of (\ref{3.1101}).\qed

\medskip

\begin{Lemma}\label{car}$($Carleson estimate$)$ Let $1<\alpha<\beta\leq2$ and let $G$ be a
$C^{1,\beta-1}$ open set with characteristics  $r_0=1$ and
$\Lambda$. Assume that $\kappa$ satisfies the  conditions in
Proposition \ref{pr2.4.}. Let     $Q=0\in\partial G$  and assume
that  $u\geq0$ is a function on $G$ which is not identical to zero,
harmonic on $G\cap B(Q,1)$ and vanishes on $  \partial G\cap B(Q,1)$
continuously. Then there exists a constant
$A_{10}=A_{10}(n,\alpha,\beta,\Lambda,C_0, C_1,C_2)$ such that
\begin{align}\label{w} u(x)\leq A_{10} u(x_0),\ \ \ \ x\in G\cap B(Q,1/2),
\end{align}
where $x_0=(0, 1/2)$ in the coordinate system $CS_Q$.
    \end{Lemma}
\noindent{\bf Proof}\  By chain arguments, we only need to prove
(\ref{w}) for $x\in G\cap B(Q,1/8)$. By multiplying a constant we
can also assume that $u(x_0)=1$. Choose $0<\gamma<\alpha/(n+\alpha)$
and define
$$B_0=G\cap B(x,2\rho(x)),\ \ \ \ B_1=B(x,\rho(x)^\gamma).$$ Set
$$\ \ \  \ B_2=B(x_0, \rho(x_0)/3),\ \ \
\ B_3=B(x_0, 2\rho(x_0)/3)$$ and
$${\tau_1}=\inf\{t>0:X_t\notin B_0\},\ \
\ \ \ {\tau_2}=\inf\{t>0:X_t\notin B_2 \}.$$

 By (\ref{3.1101}) and
scaling, we can find a constant
$\delta=\delta(n,\alpha,\beta,\Lambda,C_0,C_1,C_2) $ such that
\begin{align}\label{delta}
P_x(X_{\tau_1}\in \partial G)>\delta,\ \ \ \ x\in G\cap B(Q,1/4).
\end{align}
 By   Harnack inequality (\ref{Harnack}), there
exists $\beta'=\beta(n,\alpha,\beta,C_0,C_1,C_2) $ such that
\begin{align}\label{w12w}
u(x)< \rho(x)^{-\beta'}u(x_0),\ \ \ x\in G\cap B(Q,1/4).
\end{align}
Since $u$ is harmonic on $G\cap B(Q,1)$, we have for $x\in G\cap
B(Q,1/4)$
\begin{align}\label{ww}
u(x)=&E_x(u(X_{\tau_1})I_{X_{\tau_1} \in
B_1})+E_x(u(X_{\tau_1})I_{X_{\tau_1} \notin B_1}).
\end{align}
We  first prove  that there exists constant $l_0>0$ such that
\begin{align}\label{wsww}E_x(u(X_{\tau_1})I_{X_{\tau_1} \notin B_1})\leq
u(x_0),\ \ \ x\in G_{l_0}'\cap B(Q,1/4).\end{align}

Denote the Green function  of $(X_t)$ on an open set $U$ by $G^U$.
For $x\in G_{1/8}'\cap B(Q,1/4)$ satisfying     $$ |x-y|\leq
2|z-y|,\ \ \ \ \ z\in B_0,\ \ y\notin B_1,
$$ we have by Proposition
\ref{le5.1.8} and the last conclusion in Lemma \ref{forget}
\begin{align}\label{first} &E_x(u(X_{\tau_1})I_{X_{\tau_1} \notin
B_1})\nonumber\\=& \mathcal{A}(n,-\alpha)\int_{B_0}\int_{y\in
G,|y-x|>\rho(x)^\gamma}
\frac{\kappa(z,y)G^{B_0}(x,z)}{|z-y|^{n+\alpha}}u(y)dzdy\nonumber\\
\leq &
2^{n+\alpha}\mathcal{A}(n,-\alpha)\int_{B_0}G^{B_0}(x,z)\int_{y\in
G,|y-x|>\rho(x)^\gamma}
\frac{C_2u(y)}{|x-y|^{n+\alpha}}dzdy\nonumber\\
= & 2^{n+\alpha}\mathcal{A}(n,-\alpha)(E_x{\tau_1})\int_{y\in
G,|y-x|>\rho(x)^\gamma} \frac{C_2u(y)}{|x-y|^{n+\alpha}}dy
\nonumber\\
\leq &2^{n+2\alpha}C_2A_6'
\mathcal{A}(n,-\alpha)\rho(x)^\alpha\Big(\int_{y\in G,
|y-x|>\rho(x)^{\gamma}, |y-x_0|> 2\rho(x_0)/3}
\frac{u(y)}{|x-y|^{n+\alpha}}dy\nonumber\\+&\int_{
|y-x_0|\leq2\rho(x_0)/3}\frac{u(y)}{|x-y|^{n+\alpha}}dy\Big)
\nonumber\\
:= &
2^{n+2\alpha}C_2A_6'\mathcal{A}(n,-\alpha)\rho(x)^\alpha(I_1+I_2).
\end{align}
Similarly,
\begin{align}\label{wwws}
u(x_0)\geq &E_{x_0}(u(X_{\tau_2})I_{X_{\tau_2}\notin B_3}) \nonumber\\
\geq &
2^{-(n+\alpha)}C_1\mathcal{A}(n,-\alpha)\int_{B_2}G^{B_2}(x,z)\int_{y\in
G,|y-x|>2\rho(x_0)/3}
\frac{u(y)}{|x_0-y|^{n+\alpha}}dzdy\nonumber\\
\geq  &2^{-(n+\alpha)}C_1A_6
\mathcal{A}(n,-\alpha)(\rho(x_0)/3)^\alpha\int_{y\in
G,|y-x|>2\rho(x_0)/3} \frac{u(y)}{|x_0-y|^{n+\alpha}}dy.
 \end{align}
We have $|y-x|\geq2^{-1} \rho(x)^\gamma|y-x_0|$ if $|y-x|\geq
\rho(x)^\gamma$ and $x \in {B(Q,1/4)}$. This and (\ref{wwws}) show
that
\begin{align}\label{wwwr}I_1
\leq &2^{n+\alpha}\rho(x)^{-\gamma(n+\alpha)} \int_{y\in G,
|y-x_0|\geq 2\rho(x_0)/3}
\frac{u(y)}{|x_0-y|^{n+\alpha}}dy\nonumber\\\leq&
2^{2(n+\alpha)}3^\alpha A_6'(C_1
\mathcal{A}(n,-\alpha)\rho(x_0)^\alpha)^{-1}\rho(x)^{-\gamma(n+\alpha)}
u(x_0).\end{align} On the other hand, if $\rho(x)<\rho(x_0)/6$,  we
have by Harnack inequality (\ref{Harnack})
\begin{align}\label{wwwr}I_2
\leq &\int_{ |y-x_0|\leq2\rho(x_0)/3}\frac{u(y)}{|x-y|^{n+\alpha}}dy\nonumber\\
\leq &2^{n+\alpha}A_7\int_{ |y-x|>
\rho(x_0)/6}\frac{u(x_0)}{|x-y|^{n+\alpha}}dy
\nonumber\\
\leq &(2\pi)^n 2^{n+\alpha}A_7
(\rho(x_0)/6)^{-\alpha}u(x_0).\end{align} Combing
(\ref{first})-(\ref{wwwr}), we have for some constant
$c=c(n,\alpha,\beta,\Lambda,C_0,C_1,C_2)$
\begin{align}\label{ds}
E_x(u(X_{\tau_1})I_{X_{\tau_1} \notin B_1}) \leq& c
\rho(x)^{\alpha-\gamma(n+\alpha)}u(x_0),\ \ \ \ \ \ x\in
G_{\rho(x_0)/6}\cap B(Q,1/4).
\end{align}
Noticing that $\alpha-\gamma(n+\alpha)>0$, by choosing
$l_0=l_0(n,\alpha,\beta,\Lambda,C_0,C_1,C_2) $ small enough, we get
(\ref{wsww}) from (\ref{ds}).

 Suppose that there
exists $x_1\in G\cap B(Q,1/8)$ such that   $u(x_1)\geq M
=M(n,\alpha,\beta,\Lambda,C_0,C_1,C_2)$
$>l_0^{-\beta'}\vee(1+\delta^{-1})$ ($M$ will be fixed later). By
(\ref{w12w}), $M>l_0^{-\beta'} $ and $u(x_0)=1$ we have
$\rho(x_1)<l_0$.
  By    (\ref{ww}),(\ref{wsww}) and $M>1+\delta^{-1}$,
\begin{align}
 E_{x_1}(u(X_{\tau_1})I_{X_{\tau_1} \in B_1})\geq
 \frac{1}{1+\delta}M.\nonumber
\end{align}
From this inequality and (\ref{delta}) we   can find $x_2\in G$ such
that
$$|x_1-x_2|\leq  \rho(x_1)^{\gamma},\ \
u(x_2)>(1-\delta^2)^{-1}M.$$  Inductively, if $x_{k}\in G\cap
B(Q,1/4)$ for some  $k\geq 2$, we can find $x_{k+1}\in G$ such that
\begin{align}\label{wqw}|x_{k+1}-x_k|\leq \rho(x_k)^{\gamma},\ \
u(x_{k+1})>(1-\delta^2)^{-1}u(x_k)> (1-\delta^{2})^{-k}M.\end{align}
By (\ref{w12w}) and (\ref{wqw}), we have $\rho(x_k)\leq
(1-\delta^2)^{ k/\beta'}M^{-1/\beta' }$. Therefore, if (\ref{wqw})
holds, we have
$$|x_k|\leq |x_1|+\sum_{i=1}^{k-1}|x_{i+1}-x_i|\leq 1/8+
(1-(1-\delta^2)^{ 1/\beta'})^{-1}M^{-1/\beta'}.$$ Thus for $M=
(l_0^{-\beta'}\vee(1+\delta^{-1}))\vee(8^{\beta'}(1-(1-\delta^2)^{
1/\beta'})^{-\beta'})$, we can find $x_k\in G\cap B(Q,1/4)$
satisfying (\ref{wqw}) for all $k\geq 1$. This gives a contradiction
by noticing that $\lim_{k\rightarrow \infty}u(x_k)=\infty$ and that
$u$ vanishes on $  \partial G\cap B(Q,1)$ continuously. Therefore
$\sup_{y\in G\cap B(Q,1/8)}u(y)\leq  M$.\qed
\medskip

\begin{Remark}\label{d}
Let $G$ and   $\kappa$ satisfy  the conditions in Theorem
\ref{pr2.3.2}, then we can prove the hitting probability estimates
in Lemma \ref{pr2.3.20} and the Carleson estimate in Lemma \ref{car}
still  hold.  This is due to   that we have the same (super,sub)
harmonic functions for
$\kappa=\frac{|x-y|^{n+\alpha}}{|x-\overline{y}|^{n+\alpha}}$ and
$\kappa\equiv1$ $($see $(\ref{112}))$ and the term $C'|x-y|$ does
not destroy the (super,sub) harmonic functions which we construct
above  (c.f. Lemma \ref{pr}). We omit the proof of this extension
because it can be done by following the arguments for
 $\kappa\in C^1(\overline{{G}}\times {\overline{G}})$.  Notice that   function $\kappa$ in
  Theorem \ref{pr2.3.2} satisfies condition $[C_1,C_2,C_3,-1]$ for
some constant $C_1,C_2,C_3$, and hence
 the Harnack inequality in Theorem  \ref{pr4.1.33} holds.
\end{Remark}

Before  proving   Theorem 1.1, we give some remarks on the
assumptions of $\kappa$ and    $G$. When $\psi_2\neq0$, the
condition of $\kappa$ in $(\ref{1})$ is to  study the reflected
subordinate Brownian motion. However, due to the definition of the
reflection point, we need $C^{1,1}$ condition on $G$ in  Theorem 1.1
 when  $\psi_2\neq0$. Let $G$ be a $C^2$ open set in $\R^n$. By the
Appendix in \cite{GiT},  there exists $\delta_0>0$ such that, for
any $x\in G_{\delta_0}'$, there is an unique point $\xi(x)\in
\partial G$ satisfying  $|x-\xi(x)|=\rho(x)$,  $\xi\in
C^1(\overline{G_{\delta_0}'})$ and $\rho\in
C^2(\overline{G_{\delta_0}'})$.  For $x\in G_{\delta_0}'$, define
the reflection point   of $x$ by
\begin{align} \label{ref}\overline{x}=2\xi(x)-x.
\end{align}
When $G$ is a  $C^{1,1}$ open set,    $\xi$ and  $\rho$ are
Lipschitz and $C^{1,1}$ in a neighborhood of $\partial G$,
respectively. The proof for the uniqueness of $\xi(x)$ is similar to
\cite{GiT}. The Lipschitz and $C^{1,1}$ properties follow by
   the $C^2$ case and the standard  smooth approximation.\medskip

\textbf{ Proof of Theorem 1.1}: First we assume that $G$ is a
$C^{1,\beta-1}$ open set with characteristics  $r_0<1,\Lambda$ and
$\kappa$ satisfies the  conditions in Proposition \ref{pr2.4.}. Let
$u\geq0 $ be a function on $G$ which is not identical to zero,
harmonic on $G\cap B(Q,r)$ and vanishes continuously on $\partial
G\cap B(Q,r)$   for some $0<r<r_0$. By scaling and translation we
can assume that $r=1$ and $Q=0$. Take the  coordinate system $CS_0$
and denote
$$\ \ \ \ \ B_0=\triangle(0,A_8,A_8),\ \
 \ \ \ B_1=\triangle(0,2A_8,2A_8),\ \ \ \ \ \ \tau=\tau_{B_0}.$$ By scaling, we can
also  assume that
  \begin{align}\label{wwq}B_1\subseteq B(0,1/3).\end{align}
Write  $$x_0=(0,1/2),\ \ \ \ x_1=(0,3A_8/2).$$ By Harnack inequality
(\ref{Harnack}),  we have
\begin{align}\label{3.1101f}
k_1^{-1}u(x_0) \leq  u(x) \leq k_1 u(x_0),\ \ \ \ x\in
\triangle(0,2A_8,A_8)\setminus\triangle(0,A_8,A_8)
\end{align}
for some constant $k_1=k_1(n,\alpha,\Lambda, C_1,C_2,M)$.
 Next we assume that $ \ x\in
\triangle(0,A_8,A_8)$ with $\widetilde{x}=0$. Since $u$ is harmonic
on $G\cap B(Q,1)$,  we have by (\ref{3.1101}) and (\ref{3.1101f})
 \begin{align}\label{3.a1101}
  u(x)=&E_xu(X_{\tau})\geq  k_1^{-1}u(x_0)P_x({ u(X_{\tau})
  \in \triangle(0,2A_8,A_8) })\nonumber\\
  \geq &A_9^{-1}k_1^{-1}u(x_0)\rho(x)^{\alpha-1}.
\end{align}
By the same calculation as  (\ref{wwws}), we have
\begin{align}\label{wwwxss}
u(x_1)\geq k_2\int_{y\notin B_1} \frac{u(y)}{|x_1-y|^{n+\alpha}}dy
 \end{align}
 for some $k_2=k_2(n,\alpha,  C_1,C_2,M)$.
By   (\ref{wwq}) and Proposition \ref{le5.1.8}, we can also find a
constant
 $k_3=k_3(n,\alpha, \Lambda) $ such that (c.f. (\ref{t}))
 \begin{align}\label{use}
E_x(\tau)\leq k_3P_x(X_\tau \in G\setminus B_1).
 \end{align}
 By definition   of   $x_1$, we can find a constant
$k_4=k_4(\Lambda)$ such that $|z-y|\geq k_4 |x_1-y|$ for $z\in B_0$
and $y\notin B_1$. Thus, by Proposition \ref{le5.1.8}, Lemmas
\ref{forget}, \ref{pr2.3.20}, \ref{car} and applying
 (\ref{wwq}), (\ref{3.1101f}), (\ref{wwwxss}) and
 (\ref{use}),  we have
\begin{align}\label{firstd} u(x)=&\mathcal{A}(n,-\alpha)\int_{B_0}\int_{y\in
G\cap B_1}
\frac{\kappa(z,y)G^{B_0}(x,z)}{|z-y|^{n+\alpha}}u(y)dzdy\nonumber\\+&
\mathcal{A}(n,-\alpha)\int_{B_0}\int_{y\notin G\cap B_1}
\frac{\kappa(z,y)G^{B_0}(x,z)}{|z-y|^{n+\alpha}}u(y)dzdy\nonumber\\
\leq & C_2\mathcal{A}(n,-\alpha)A_{10}u(x_0)P_x(X_\tau\in B_1)
+C_2\mathcal{A}(n,-\alpha)k_4^{-(n+\alpha)}\int_{B_0}\int_{y\notin
B_1}
\frac{ G^{B_0}(x,z)}{|x_1-y|^{n+\alpha}}u(y)dzdy\nonumber\\
\leq  & C_2\mathcal{A}(n,-\alpha)A_9A_{10}u(x_0)\rho(x)^{\alpha-1}
+C_2\mathcal{A}(n,-\alpha) k_4^{-(n+\alpha)}E_x(\tau) \int_{y\notin
B_1} \frac{u(y)}{|x_1-y|^{n+\alpha}}dy
\nonumber\\
\leq  &
C_2\mathcal{A}(n,-\alpha)u(x_0)A_9(A_{10}+k_1k_2^{-1}k_3k_4^{-(n+\alpha)})\rho(x)^{\alpha-1}.
\end{align}
Combing (\ref{3.a1101}) and (\ref{firstd}) we prove  (\ref{yesr3})
for $ x\in \triangle(0,A_8,A_8)$ with  $\widetilde{x}=0$. By
considering the  coordinate system   $CS_y$ for $y\in
B(Q,2/3)\cap\partial G$, applying the arguments above and  the
Harnack inequality (\ref{Harnack}), we  can prove (\ref{yesr3}). The
general case  can be proved similarly with the help of Remark
\ref{d}.

\qed\medskip

\section{\normalsize  Boundary Harnack inequality of
$\Delta^{{\alpha}/{2},\kappa}_G$ on Lipschitz domain}

To simplify notations, we assume that $\kappa\equiv1$ in the
arguments below because the estimates are the same for the general
cases.
 Let $G$ be a Lipschitz domain with
characteristic $r_0$ and $\Lambda$, i.e.,  for each $x_0\in
\partial
 G$, we can find a Lipschitz function $\Gamma_{x_0}:\R^{n-1}\rightarrow
 \R$ with Lipschitz coefficient not greater than $\Lambda$ and an
 orthonormal coordinate system $CS_{x_0}$ with which it holds that
 \begin{align}
G\cap B \bigr(x_0,r_0\bigr)=\{y=(y_1,\cdots,y_n):\ y_n>
\Gamma_{x_0}(y_1,\cdots y_{n-1}) \ \}\cap B\bigr( x_0,r_0\bigr).
\end{align}
The following hitting probability estimate  is obvious  for the
Brownian motion case. Here we use  capacity to give this estimate.
We refer to \cite{FOT} for more details about  capacity and  energy
measure class $S_{00}$ of   symmetric Markov processes.

\begin{Lemma}\label{carqwrwer2}
Let $G$ be a Lipschitz domain with characteristic $r_0=1$ and
$\Lambda$. Let $x_0=0\in \partial G$ and choose the coordinate
system $CS_0$.  Assume that $A$ is  a constant such that
$\triangle(0,A,A)\subset G\cap B(0,1)$. Then there exists a constant
$A_{11}=A_{11}(n,\alpha,\Lambda)$ such that
\begin{align}\label{lfffl}
P_x( X_{\tau}\in \partial G)\geq  1/A_{11},\ \ \ \
x=(\widetilde{0},x_n),\ \ 0<x_n<A/2,
\end{align}
where $\tau=\inf\{t>0: X_t\in \triangle(0,A,A)^c\}$.
    \end{Lemma}
\noindent{\bf Proof}\  By scaling, we may assume that
$A>1/(3(1+\Lambda))$ without loss of generality.  Denote the
  heat
kernel of $(X_t)$ by $p(t,x,y)$. By Theorem 1.1 in \cite{CHKB},
there exists constant $k_1=k_1(n,\alpha,\Lambda)$ such that
\begin{align}\label{asa}
k_1\big(t^{-n/\alpha}\wedge \frac{t}{|x-y|^{n+\alpha}}\big)\leq
p(t,x,y) \leq k_1^{-1}\big( t^{-n/\alpha}\wedge
\frac{t}{|x-y|^{n+\alpha}}\big),\ \ \ 0<t<1.
\end{align}
Set $F= \overline{{\triangle(0,A,A)}^c\cap G}$ and denote by $(Y_t)$
the killed process of $(X_t)$   when hitting $F$. We know that the
heat kernel  $p^0(t,x,y)$  of $(Y_t)$ is given by
\begin{align}\label{ii}
p^0(t,x,y)=p(t,x,y)-\int_0^t\int_Fp(t-s,z,y)P_x(X_{\sigma}\in dz,
\sigma\in ds),
\end{align}
where $\sigma=\inf\{t>0:X_t\in F\}$. Noticing that
$A>1/(3(1+\Lambda))$ and  choosing $\delta$ small enough we get by
(\ref{asa}) and (\ref{ii})
\begin{align}\label{ppm}
p^0(t,x,y)\geq k_2\big(t^{-n/\alpha}\wedge
 \frac{t}{|x-y|^{n+\alpha}}\big),\ \ \
0<t<1,\ x,y\in \overline{G}\cap B(0,\delta)
\end{align}
for some $k_2=k_2(n,\alpha,\Lambda)$. Set $\Gamma=\partial G \cap
\overline{B(0,\delta)}$.
 Define the
1-potential kernel of $(Y_t)$ by $U_1^0(x,y)=\int_0^\infty
e^{-t}p^0(t,x,y)\ dt$ and define for measure $\mu$ on
$\overline{G}\setminus F$ \begin{align}\label{hitting}U_1^0
\mu(x)=\int_{\overline{G}\setminus F}U_1^0(x,y)\mu(dy).
\end{align}
 By Theorem 4.2.5 in \cite{FOT} and the continuity argument,  there exists a 1-equilibrium measure
$\nu_\Gamma$ supported on $\Gamma$ such that
\begin{align}\label{sb}
U_1^0\nu_\Gamma(x)=E_x^0(e^{-\sigma_\Gamma}), \ \ \ \ \  \ \ \ x\in
\overline{G}\setminus F,
\end{align}
where $\sigma_\Gamma=\inf\{t>0: Y_t\in \Gamma \}$. By problem 2.2.2
in \cite{FOT},
\begin{align} \nu_\Gamma(\Gamma)=\sup\{\mu(\Gamma): \mu\in S_{00},\
supp[\mu]\subseteq K,\ U_1^0\mu \leq 1\}.
\end{align}
Direct calculations shows that $\mu=\delta I_{\Gamma}m(dx)\in
S_{00}$ for $\delta>0$. Choosing $\delta$ small enough and applying
the second inequality in (\ref{asa}), we get
$\nu_\Gamma(\Gamma)>k_3(n,\alpha,\Lambda)$. Therefore by (\ref{ppm})
and (\ref{sb}), we have for $x \in  \overline{G}\cap B(0,\delta)$
\begin{align}\label{sbbf}
E_x^0(e^{-\sigma_\Gamma})\geq & e^{-1}\int_{\Gamma}\int_0^1p_0(t,x,y)dt\nu_\Gamma(dy)\nonumber\\
\geq & e^{-1}\nu_\Gamma(\Gamma)\inf_{y\in
\Gamma}\int_0^1p_0(t,x,y)dt\geq k_4(n,\alpha,\Lambda),
\end{align}
which implies that
\begin{align}\label{sbfg}
P_x^0(\sigma_\Gamma<\infty) \geq E_x^0(e^{-\sigma_\Gamma})\geq k_4.
\end{align}
Noticing that $P_x(X_\tau\in \partial G)\geq
P_x^0(\sigma_\Gamma<\infty) $, we get (\ref{lfffl}) for $x \in
\overline{G}\cap B(0,\delta)$. Thus we complete the proof by the
Harnack inequality in Theorem \ref{pr4.1.33}. \qed
\medskip

 Let  $A$ be a constant such that
$\triangle(0,2A,2A)\subset G\cap B(0,1)$ under the coordinate system
$CS_0$ for  a Lipschitz domain $G$ with characteristic $r_0=1$ and
$\Lambda$.  Set
\begin{align}\label{A} K_0=&\triangle(0,A,A),\ \ \
K_1=\triangle(0,2A,A)\setminus K_0,\ \ \ K_2=G\setminus \triangle(0,2A,A);\\
\label{AA}
 H_1=&\{X_{\tau_{K_0}}\in K_1\},\ \ \ \ \ H_2=\{X_{\tau_{K_0}}\in
 K_2\}.
\end{align}

\begin{Lemma}\label{car54}  With notations
defined  in (\ref{A}) and (\ref{AA}), for any   $k\geq 0$, there
exists a constant $A_{12}=A_{12}(n,\alpha, \Lambda, A,k)$ such that
\begin{align}\label{ll}
P_y( H_1)\geq A_{12}\rho(y)^\alpha |\ln \rho(y)|^k,\ \ \ \ y\in K_0,
\ |\widetilde{y}|\leq A/2.
\end{align}
    \end{Lemma}
\noindent{\bf Proof}\  Let $y\in K_0$ with $|\widetilde{y}|\leq A/2$
and $\tau=\tau_{B(y,\rho(y)/2)}$. We   assume also that
$B(y,\rho(y)/2)\subset K_0$. Otherwise  (\ref{ll}) can  be verified
by showing that $P_y(H_1)>k_1(n,\alpha,\Lambda,A)$. As the
calculations in (\ref{t}), we have by Lemma \ref{forget}
\begin{align}\label{tq} P_y(H_1)
\geq &  k_2(n,\alpha)E_{y}({\tau})\int_{z\in K_1 }
\frac{1}{|z-y|^{n+\alpha}}dz
\nonumber\\
\geq &kk_3(n,\alpha,\Lambda,A)\rho(y)^\alpha, \end{align} which
gives (\ref{ll}) for $k=0$. Suppose (\ref{ll}) holds for some $
k\geq 0$. By the strong Markov property
\begin{align}\label{dsfs} P_y(H_1)
\geq & P_y( \tau<\tau_{K_0}  , \ X_{  \tau_{K_0}  }  \in
K_1)\nonumber\\
=& \int_{K_0}P_z(H_1) P_y(X_\tau\in dz) \nonumber\\\geq&
 k_4(n,\alpha)\int_{z\in K_0\setminus B(y,\rho(y)/2)}P_z(H_1)
\frac{E_{y}({\tau})}{|z-y|^{n+\alpha}}dz \nonumber\\\nonumber\\\geq&
 k_5(n,\alpha,\Lambda,A,k)\int_{z\in K_0\setminus B(y,\rho(y)/2), |\widetilde{z}|\leq A/2}
\frac{\rho(y)^{\alpha}\rho(z)^{\alpha} |\ln
\rho(z)|^k}{|z-y|^{n+\alpha}}dz.
\end{align} Direct calculation shows
that $ P_y(H_1)\geq A_{12}'\rho(y)^\alpha |\ln \rho(y)|^{k+1}$.
Hence the proof completes  by induction. \qed \medskip

\begin{Lemma}\label{caradfaadf2123}   With notations
defined  in (\ref{A}) and (\ref{AA}), there exists a constant
$A_{13}=A_{13}(n,\alpha,\Lambda)$ such  that
\begin{align}\label{qq}
P_y( H_2)\leq A_{13}P_y( H_1),\ \ \ \ y\in K_0, \ \widetilde{y}=0.
\end{align}
    \end{Lemma}
\noindent{\bf Proof}\ To simplify the arguments, we assume that $G$
is a special Lipschitz domain. By scaling, we assume that $A=1$ in
(\ref{A}). For $i\geq 1$, set
$$J_i=\triangle(0,2^{-i},r_i)\setminus  \triangle(0,2^{-i-1},r_i),\ \ \ \ r_i=\frac{1}{2}-\frac{1}{50}\sum_{j=1}^i\frac{1}{j^2},$$
and $r_0=r_1$.  Define for $i\geq 1$ \begin{align}\label{q1q}
d_i=\sup_{z\in J_i}P_z( H_2)/P_z( H_1),\ \ \
\widetilde{J}_i=\triangle(0,2^{-i},r_{i-1}),\ \ \ \
\tau_i=\tau_{\widetilde{J}_i}.
\end{align}
By Harnack inequality, each $d_i$ is finite.  Noticing that
$\tau_i\leq \tau_{K_0} $ and applying  the strong Markov property,
we have for $z\in J_i$ and $i\geq 2$
\begin{align}\label{q12q}  P_z( H_2)=& P_z(X_{\tau_{K_0}}\in K_2, \ X_{\tau_i}\in \cup_{k=1}^{i-1}J_k)
+P_z(X_{\tau_{K_0}}\in K_2, \ X_{\tau_i}\in G\setminus
\cup_{k=1}^{i-1}J_k)\nonumber\\
\leq  &\sum_{k=1}^{i-1} \int_{J_k} P_z (X_{\tau_i}\in dw) P_{w}(H_2)
+P_z( \ X_{\tau_i}\in G\setminus \cup_{k=1}^{i-1}J_k)
\nonumber\\
\leq  &\sum_{k=1}^{i-1}  d_k\int_{J_k} P_z (X_{\tau_i}\in dw)
P_{w}(H_1) +P_z(  \ X_{\tau_i}\in G\setminus \cup_{k=1}^{i-1}J_k)
\nonumber\\
\leq  &(\sup_{1\leq k\leq i-1}  d_k)  P_{z}(H_1) +P_z( X_{\tau_i}\in
G\setminus \cup_{k=1}^{i-1}J_k).
\end{align}
Define $\sigma_0=0, \sigma_1=\inf\{t>0: |X_t-X_0|\geq 2^{-i}\} $ and
define by induction $\sigma_{m+1}=\sigma_1\circ\theta_{\sigma_{m}}$
for $m\geq 1$. Similar to the calculation of (\ref{t}), we have for
some constant $k_1<1$ independent of $i$ such that
$$P_{w}(X_{\sigma_1}\in \widetilde{J}_i)\leq 1-P_{w}(X_{\sigma_1}\in  \cup_{k=1}^{i-1}J_k)<k_1,\ \ \ w\in  \widetilde{J}_i.
$$ Therefore, for   $z\in J_i$ and positive integer $l$, we have  by the strong Markov property for
\begin{align}\label{q132q}    P_z( \tau_{i}>\sigma_{li})\leq& P_z(X_{\sigma_{k}}\in \widetilde{J}_i, 1\leq k\leq
li
 )\nonumber\\=&\int_{w\in  \widetilde{J}_i}P_z(X_{\sigma_{k}}\in \widetilde{J}_i, 1\leq k\leq
li-2, X_{\sigma_{li-1}}\in dw
 ) P_w(X_{\sigma_1}\in  \widetilde{J}_i)
 \nonumber\\      \leq& P_z(X_{\sigma_{k}}\in \widetilde{J}_i, 1\leq k\leq
li-1
 )k_1\leq k_1^{li}.
\end{align}
Recall that $\widetilde{x}$ is the first $n-1$ coordinate of $x$. On
$\{X_{\tau_{i}}\in G\setminus \cup_{k=1}^{i-1}J_k,\ \ \tau_{i}\leq
\sigma_{li}\}$ with $X_0=z\in J_i$,  we have
$|\widetilde{X}_{\sigma_k}-\widetilde{X}_{\sigma_0}|>
\frac{1}{50i^2}-2^{-i}$ for some $1\leq k\leq li$ which implies for
some $1\leq k'\leq li$ $$ |X_{\sigma_{k'}}-X_{\sigma_{k'-1}}|\geq
(\frac{1}{50i^2}-2^{-i})/(li).
$$
Therefore, we have for some $i_0\geq 2$
\begin{align}& \{X_{\tau_{i}}\in G\setminus
\cup_{k=1}^{i-1}J_k,\ \ \tau_{i}\leq \sigma_{li}\}\nonumber\\
 \subseteq &   \cup_{k=1}^{li}\{|
{X}_{\sigma_k}- {X}_{\sigma_{k-1}}|\geq 1/(100li^3),
{X}_{\sigma_{k-1}}\in \widetilde{J}_{i}
 \},\ \ \ i\geq i_0,
\end{align}
and hence
 \begin{align}\label{dsa}& P_z(X_{\tau_{i}}\in G\setminus
\cup_{k=1}^{i-1}J_k,\ \ \tau_{i}\leq \sigma_{li})\nonumber\\
 \leq  &\sum_{k=1}^{li} P_z(|
{X}_{\sigma_k}- {X}_{\sigma_{k-1}}|\geq 1/(100li^3),
{X}_{\sigma_{k-1}}\in \widetilde{J}_{i}
 )
 \nonumber\\
 \leq  &li \sup_{z\in \widetilde{J}_{i}} P_z(|
{X}_{\sigma_1} |\geq 1/(100li^3))
\nonumber\\
 \leq  &k_2li  2^{-\alpha i}(li^{3})^{ \alpha}.
\end{align}
The proof of the last inequality above is similar to (\ref{first}).
By (\ref{q132q}), (\ref{dsa}) and  choosing $l$ big enough, we have
for $z\in J_i$ and $i\geq i_0$,
\begin{align}\label{as}
P_z( X_{\tau_{i}}\in G\setminus \cup_{k=1}^{i-1}J_k ) \leq
k_1^{li}+k_2li  2^{-\alpha i}(li^{3})^{ \alpha}\leq k_3   2^{-\alpha
i}i^{3 \alpha+1}. \end{align} By (\ref{q12q}), (\ref{as}) and Lemma
\ref{car54}, for $z\in J_i$ and $i\geq i_0$
\begin{align}  P_z( H_2)/P_z(H_1)& \leq  \sup_{1\leq k\leq i-1}  d_k   +P_z( X_{\tau_i}\in
G\setminus \cup_{k=1}^{i-1}J_k)/P_z(H_1) \leq \sup_{1\leq k\leq i-1}
d_k+k_4/i^2.\nonumber
\end{align}
This implies that
\begin{align}  d_i& \leq  \sup_{1\leq k\leq i_0-1} d_k+ k_4\sum_{k=1}^i1/k^2\leq
\sup_{1\leq k\leq i_0-1} d_k+3k_4,\nonumber
\end{align}
which completes the  proof of this lemma. \qed \medskip
\begin{Remark}
One may use the method in \cite{BOG} to give a better estimate of
(\ref{ll}). The proof of Lemma \ref{caradfaadf2123} is an adaption
of   the Brownian motion case.
\end{Remark}

\textbf{ Proof of Theorem 1.2}: In the proof of Lemma \ref{car}, we
only use the $C^{1,\beta-1}$ property in (\ref{delta}). Thus   we
can prove the Carleson estimate for the Lipschitz case with Lemma
\ref{carqwrwer2} in place of (\ref{delta}). Therefore,  we can prove
Theorem 1.2 by the standard arguments of BHI with the help of
Theorem \ref{pr4.1.33} and Lemma \ref{caradfaadf2123}. \qed

\section{\normalsize  Boundary Harnack inequality of $\Delta^{{\alpha}/{2}}$}

When $G=\R^n$,  $\Delta^{{\alpha}/{2}}_{G}$ is the fractional
Laplacian $\Delta^{{\alpha}/{2}}$.  Recall that $w_p(y)=y_n^{ p}$
for $ y\in\R^n_+$.  We extend   these functions to $\R^n $ by taking
zero on $  \R^n\setminus \R^n_+ $.
 Next we give a formula of  $\Delta^{ {\alpha}/{2}}w_p$. Integration by parts
formula shows that for $0<p<\alpha$
\begin{align}\label{654}
&\int_0^\infty\frac{y^{p}-x^{p}}{|y-x|^{1+\alpha}}dy
=x^{p-\alpha}\int_0^\infty\frac{y^{p}-1}{|y-1|^{1+\alpha}}dy\nonumber
\\
=&\lim_{\varepsilon\downarrow
0}x^{p-\alpha}\left(\int_0^{1-\varepsilon}\frac{y^{p}-1}{|y-1|^{1+\alpha}}dy
+\int_{1+\varepsilon}^\infty\frac{y^{p}-1}{|y-1|^{1+\alpha}}dy\right)
\nonumber
\\
=&\frac{1}{\alpha}x^{p-\alpha}+\frac{p}{\alpha}
x^{p-\alpha}\int_0^{1}\frac{y^{\alpha-p-1}-y^{p-1}}{|y-1|^{\alpha}}dy.
\end{align} Thus, for $n=1$ we
have
\begin{align}
\Delta^{\alpha/2}w_{p}(x)&=\mathcal{A}(1,-\alpha) \frac{p}{\alpha}
x^{p-\alpha}\int_0^{1}\frac{y^{\alpha-p-1}-y^{p-1}}{|y-1|^{\alpha}}dy,\
\ 0<p<\alpha .
\end{align}
Applying spherical coordinate transform from $(y_1,\cdots,y_n)\in
\R^n$ to $(r,\theta_1, \cdots,\theta_{n-1})\in
[0,\infty)\times[0,\pi]^{n-2}  \times [0,2\pi)$, this gives for  $n>
1$
\begin{align}\label{4567}
&\mathcal{A}(n,-\alpha)^{-1}\Delta^{\alpha/2}w_{p}(x)\nonumber\\=&\lim_{\varepsilon\downarrow
0} \int_{|y-x|>\varepsilon} \frac{w_p(y)-w_p(x)}{|x-y|^{n+\alpha}}
dy\nonumber\\=&\lim_{\varepsilon\downarrow 0}
             \int_0^{\pi/2}
 d\theta_1\cdots
 \int_0^\pi
 d\theta_{n-2}\int_0^{2\pi}
 \varphi(\theta_1,\cdots\theta_{n-2})\ d\theta_{n-1}\cdot\nonumber\\
 &\int_{-\infty}^\infty
I_{\{|t-\frac{x_n}{\cos\theta}|>\varepsilon\}}
             \cos^{p}\theta_1\frac{t^{p}I_{t\geq0}-(\frac{x_n}{\cos \theta_1})^{p}}
               {|t-\frac{x_n}{\cos\theta_1}|^{1+\alpha}} \
               dt
               \nonumber\\=&
               \frac{p}{\alpha }\int_0^{1}\frac{y^{\alpha-p-1}-y^{p-1}}{|y-1|^{\alpha}}dy
             \int_0^{\pi/2}
 d\theta_1\cdots
 \int_0^\pi
 d\theta_{n-2}\int_0^{2\pi}
 \varphi(\theta_1,\cdots\theta_{n-2})\
             (\cos^{\alpha}\theta_1)
x^{p-\alpha}_n d\theta_{n-1}
               \nonumber\\
               =&\frac{p}{\alpha }\int_0^{1}\frac{y^{\alpha-p-1}-y^{p-1}}{|y-1|^{\alpha}}dy
 \int_{|y|=1, y_n\geq 0}y_n^{\alpha}\  m(dy)\cdot w_{p-\alpha}(x),
\end{align}
where  $
 \varphi(\theta_1,\cdots\theta_{n-2})=\sin^{n-2}\theta_1
\sin^{n-3}\theta_2\cdots\sin\theta_{n-2}$,   $m(dy)$ is the
$(n-1)$-dimensional Lebesgue measure and we use the following
transform in the calculation above $$r=t-x_n/\cos\theta_1,\ \
\theta_1\in[0,\pi/2); \ \ \ \ \ \ \ -r=t+x_n/\cos\theta_1,\ \
\theta_1\in(\pi/2,\pi].
$$  Denote for $0<p<\alpha$ and $n\geq 1$
\begin{align}\Lambda(n,\alpha,p)=&
\frac{p\mathcal{A}(n,-\alpha)}{\alpha }
\int_0^{1}\frac{y^{\alpha-p-1}-y^{p-1}}{|y-1|^{\alpha}}dy
 \int_{|y|=1,y_n\geq 0}y_n^{\alpha}\  m(dy),\nonumber  \\
 \overline{\Lambda}(n,\alpha,p)=&
\frac{ \mathcal{A}(n,-\alpha)}{\alpha
}\left(1+p\int_0^{1}\frac{y^{\alpha-p-1}-y^{p-1}}
{|y-1|^{\alpha}}dy\right)
 \int_{|y|=1,y_n\geq 0}y_n^{\alpha} \ m(dy)\nonumber\end{align}
 with convention that $ m(dy) $ is the Dirac measure  for $n=1$.
By   (\ref{4567}), we have the following Lemma. We notice that the
case $p=\alpha/2$ below has been obtained  in Example 3.2 of
 \cite{BAB}.
\begin{Lemma}\label{this}
Let $0<p<\alpha<2$, we have
\begin{align}\label{lag}
\Delta^{\alpha/2} w_p=&\Lambda(n,\alpha,p)w_{p-\alpha},\ \ \ \
x\in\R^n_+,
\\\label{vuyd}
\Delta^{\alpha/2}_{\R^n_+}w_p=&\overline{\Lambda}(n,\alpha,p)w_{p-\alpha},\
\ \ \ x\in\R^n_+.
\end{align}
\end{Lemma}
 Formula  (\ref{vuyd}) is another version of (\ref{112}) for $0<p<\alpha$.
 By   Lemma \ref{this} and (\ref{112}) we see that
\begin{align} \label{PPPPP}\Delta^{\alpha/2} w_p<0,\ \ -1<p<\alpha/2;\ \
\Delta^{\alpha/2} w_p=0,\ \ p=\alpha/2;\ \ \Delta^{\alpha/2} w_p>0,\
\ \alpha/2<p<\alpha.
\end{align}
Let $\kappa$ be a symmetric function  on $\R^n\times\R^n$
   such that
\begin{align}\label{dprf}
R_1<\kappa(x,y)<R_2,\ \ \ |\kappa(x,y)-\kappa(x,x)|\leq R_3|x-y|,\ \
\ x,y\in \R^n
\end{align}
for some constants $R_1,R_2,R_3>0$. In what follows, notation
$(X_t)$ is for the stable-like process on $\R^n$ associated with
$\Delta^{\frac{\alpha}{2},\kappa}$. Harmonic functions of $(X_t)$ is
again defined by (\ref{hc}).
\begin{Lemma}\label{pr2.3dfsasd.2} Let $0<\alpha\leq 1\vee\alpha<\beta\leq2$ and
   $D$
a $C^{1,\beta-1}$ open set in $\R^n$ with characteristics $r_0=1$
and  $\Lambda$. Let $\kappa$ be a symmetric function on
$\R^n\times\R^n$ satisfying  (\ref{dprf}). Then for $\alpha/2\leq
p<\alpha$ and $Q\in \partial D$, there exist function $u_p$ and
constants $A_{13}=A_{13}(\Lambda)$,
$A_{14}=A_{14}(n,\alpha,\beta,\Lambda,p,R_1,R_2,R_3)$ such that
\begin{align}\label{r*}A_{13}^{-1}I_{D\cap B(Q,2/3)}\rho(x)^p\leq u_p(x)\leq A_{13}I_{D\cap
B(Q,2/3)}\rho(x)^p,\ \ \ x\in \R^n,\end{align}
\begin{align}\label{2df1111} \Delta^{\frac{\alpha}{2},\kappa} u_p (x)
\geq A_{14} \rho(x)^{p-\alpha},\  \ \ \   x \in D \cap
B(Q,1/A_{14}),\ \ \alpha/2<p<\alpha,
\end{align}
and
\begin{align}\label{no1fd1}
|\Delta^{\frac{\alpha}{2},\kappa} u_{\alpha/2} (x)|\leq\left\{
\begin{array}{l@{\quad \quad}l}  A_{14}
\rho(x)^{\beta-\alpha/2-1},\  \ \   x \in D \cap B(Q,1/2),\ \ \
  \beta<1+\alpha/2,\\
  A_{14}|\log \rho(x)| ,\  \ \  \ \ \   x \in D \cap
B(Q,1/2),\ \ \
  \beta=1+\alpha/2,\\
A_{14} ,\  \ \ \ \ \ \ \ \ \ \ \ \ \ \ \ \   x \in D \cap B(Q,1/2),\
\ \   \beta> 1+\alpha/2.
        \end{array}                  \right. \end{align}
\end{Lemma}
\noindent{\bf Proof}\ Following the calculations in Lemma
\ref{hpr2.3.2}, Lemma \ref{pr2.3.2y} and Proposition \ref{pr2.4.} we
can prove this lemma with the help of (\ref{lag}) and (\ref{PPPPP}).
We omit the details of the proof  because, by noticing that $w_p=0$
on $\R^n\setminus\R^n_+$,  the calculation is essentially on $D$
which is the same as the regional fractional  Laplacian case.
 \qed
\medskip

 By Lemma \ref{pr2.3dfsasd.2} and  following the arguments
in Proposition \ref{pr2.3.20} and Theorem 1.1, we can prove the
following results.
\begin{Lemma}
\label{pr2dsfx.3.20dasx} Let $0<\alpha\leq 1\vee\alpha<\beta\leq2$
and    $D$   a $C^{1,\beta-1}$ open set in $\R^n$ with
characteristics $r_0=1$ and $\Lambda$. Assume that $Q=0\in \partial
D$ and $\kappa$ is a symmetric function on $\R^n\times\R^n$
satisfying (\ref{dprf}). Then there exist constants
$A_{15}=A_{15}(n,\alpha,\beta,\Lambda,R_1,R_2,R_3)<1/2$ and
$A_{16}=A_{16}(n,\alpha,\beta,\Lambda, R_1,R_2,R_3)$ such that
\begin{align}\label{3.110asd1}
 A_{16}^{-1} \rho(x)^{\alpha/2} \leq &P_x\{X_{\tau_{\triangle(0,A_{15},A_{15})}}\in
\triangle(0,2A_{15},A_{15})\}\nonumber\\\leq&
P_x\{X_{\tau_{\triangle(0,A_{15},A_{15})}}\in D\}\ \leq A_{16}
\rho(x)^{\alpha/2}
\end{align}   for $ x\in \triangle(Q,A_{15},A_{15})$ with
$\widetilde{x}=0$ under $\mbox{CS}_Q$.
\end{Lemma}

\begin{theorem}\label{guan} Assume that  $\alpha$, $\beta$, $D$ and $\kappa$
satisfy the same conditions as in Lemma \ref{pr2dsfx.3.20dasx}.  Let
$Q\in\partial D$ and $r\in (0,r_0)$. Assume that
 $u\geq0$ is a function on $D$ which is not identical to
zero, harmonic on $D\cap B(Q,r)$ and vanishes   on $  D^c\cap
B(Q,r)$. Then there exists    constant
$C=C(n,\alpha,\Lambda,R_1,R_2,R_3)$ such that
\begin{align}\label{yesr} \frac{u(x)}{u(y)}\leq
C\frac{\rho(x)^{\alpha/2}}{\rho(y)^{\alpha/2}},\ \ \ for\ x,y\in
D\cap B(Q,r/2).
\end{align}
\end{theorem}

\begin{Remark}
By taking $G=D$,    all  the conclusions in Section 3 can be
extended to $\Delta^{\frac{\alpha}{2},\kappa}$ in a similar  way,
where the reflected stable process is replaced by the stable-like
process. The Carleson estimate for
$\Delta^{\frac{\alpha}{2},\kappa}$ can be proved by the same method
as in Lemma \ref{car}. We remark that to prove the boundary Harnack
principle of $\Delta^{\frac{\alpha}{2},\kappa}$ on open sets, we
need the  method in \cite{BOGA} to get the Carleson estimate, where
the explicit Poisson kernel can be replaced by the sharp estimates
as in \cite{CS}. Theorem \ref{guan} can be generalized to operator
$\Delta^{\frac{\alpha}{2},\kappa}_G$ when we further assume that
$D\subset \overline{D}\subset G$. The proof of this generalization
is the same as the case $G=\R^n$ except that the constant depends
also on the distance between $D$ and $\partial G$.
\end{Remark}

\begin{Remark}
Since $((w_{p})_{p<1}$, $(w_p)_{p>1})$ $w_1$ is the (super,
sub)harmonic function of Laplacian on half spaces, by the Harnack
inequality in \cite{SV} and the method in this paper, we can prove
the  explicit BHI for $\Delta+\Delta^{\alpha/2}$ on  $C^{1,1}$ open
sets which gives $\rho(x)$ order decay for harmonic functions near
the boundary. With the help of this fact  we can prove that the
Green functions   of $\Delta+\Delta^{\alpha/2}$ and $\Delta$ are
comparable on a $C^{1,1}$ bounded open set.
\end{Remark} \noindent$\mathbf{ Acknowledgement}$\ \ This
paper was reported partly in the second international conference on
stochastic analysis and its applications in Seoul. The last version
of this paper without the Lipschitz case was completed when the
author was at the school of mathematics of Loughborough University.

\end{document}